\newif\ifpreprint
\preprinttrue  

\ifpreprint

\documentclass[11pt]{article}
\usepackage{fullpage}
\title{Learning Algorithm Hyperparameters \\for Fast Parametric Convex Optimization}

\author{Rajiv Sambharya$^1$ and Bartolomeo Stellato$^2$}
\date{%
    $^1$University of Pennsylvania\\
    $^2$Princeton University\\[1em]%
    \today
}

\usepackage{amsmath,enumitem,mathtools,amsthm,amssymb}

\renewcommand\arraystretch{1.2}

\newtheorem{theorem}{Theorem}

\usepackage[colorlinks,
            linkcolor=cyan,
            citecolor=cyan,
            urlcolor=cyan,
            bookmarks]{hyperref}
\let\citep\cite
\let\citet\cite

\else 

\documentclass[review,onefignum,onetabnum]{siamonline220329}

\usepackage{lipsum}
\usepackage{amsfonts}
\usepackage{graphicx}
\usepackage{epstopdf}
\usepackage{bibentry}
\ifpdf
  \DeclareGraphicsExtensions{.eps,.pdf,.png,.jpg}
\else
  \DeclareGraphicsExtensions{.eps}
\fi

\usepackage{enumitem}
\setlist[enumerate]{leftmargin=.5in}
\setlist[itemize]{leftmargin=.5in}

\newsiamremark{remark}{Remark}
\newsiamremark{hypothesis}{Hypothesis}
\crefname{hypothesis}{Hypothesis}{Hypotheses}
\newsiamthm{claim}{Claim}

\headers{Learning for Fast Parametric Convex Optimization}{R. Sambharya, and B. Stellato}

\title{Learning Algorithm Hyperparameters for Fast Parametric Convex Optimization\thanks{Submitted to the editors 11/23/2024.}}

\author{Rajiv Sambharya\thanks{University of Pennsylvania, Philadelphia, PA (\email{sambhar9@seas.upenn.edu})} \and Bartolomeo Stellato\thanks{Princeton University, Princeton, NJ (\email{bstellato@princeton.edu})}.}

\usepackage{amsopn}

\let\citep\cite
\let\citet\cite

\nocite{bertsekas1999nonlinear}
\nocite{burton2010history}
\nocite{callebaut}
\nocite{NoceWrig06}

\renewcommand\arraystretch{1.2}

\usepackage{caption}
\usepackage{graphicx}
\usepackage{subcaption}
\usepackage{hyperref}
\fi
\usepackage{tikz}
\usepackage{verbatim}
\usepackage{multicol}
\usepackage{algorithm}
\usepackage[noend]{algpseudocode} 

\usepackage{subcaption}
\usepackage{makecell}
\usepackage{csvsimple}	
\usepackage{booktabs}   
\usepackage[flushleft]{threeparttable}  
\usepackage[round-mode=places, round-integer-to-decimal, round-precision=4,
    table-format = 1.4,
    table-number-alignment=center,
    round-integer-to-decimal]{siunitx}
\usepackage{adjustbox}  

\newcommand{\bnote}[1]{}
\renewcommand{\bnote}[1]{\textcolor{red}{\textbf{B. #1}}}

\usepackage{etoolbox}

\newcommand{\fplen}{n}
\newcommand{\scsnprimal}{q}

\DeclareMathOperator*{\argmin}{argmin}

\newcommand*{\startlegend}{-0.2}
\newcommand*{\enlegend}{0.7}
\newcommand{\risk}{r_\mathcal{X}}
\newcommand{\emprisk}{\hat{r}_{S}}

\DeclareDocumentCommand{\T}{ O{z} O{} }{T\IfValueT{#2}{(#1,\theta_{#2})}\IfNoValueT{#2}{(#1,\theta)}}
\DeclareDocumentCommand{\Tj}{ O{k} O{z} O{} }{T^{#1}\IfValueT{#3}_{\theta_{#3}}{(#2)}}

\DeclareDocumentCommand{\CB}{ O{} }{C_{B_{\IfValueTF{#1}{\theta_{#1}}{\theta}}}}
\DeclareDocumentCommand{\RB}{ O{} }{R_{B_{\IfValueTF{#1}{\theta_{#1}}{\theta}}}}

\newcommand{\eg}{{\it e.g.}}
\newcommand{\ie}{{\it i.e.}}

\newcommand{\reals}{{\mbox{\bf R}}}

\newcommand{\symm}{{\mbox{\bf S}}}  

\newcommand{\Tr}{\mathop{\bf tr}}
\newcommand{\diag}{\mathop{\bf diag}}

\newcommand{\prox}{\textbf{prox}}

\newcommand{\Sec}{Section}
\newcommand{\Thm}{Theorem}

\newcommand{\Eqn}{Equation}

\newcommand{\cblock}[3]{
  \hspace{-1.5mm}
  \begin{tikzpicture}
    [
    node/.style={square, minimum size=10mm, thick, line width=0pt},
    ]
    \node[fill={rgb,255:red,#1;green,#2;blue,#3}] () [] {};
  \end{tikzpicture}%
}

\newcommand{\ccircle}[3]{%
  \raisebox{0.55\height}{
  \begin{tikzpicture}[baseline=(node.base)]%
    \node[circle, fill={rgb,255:red,#1;green,#2;blue,#3}, yshift=14] (node) at (0,1) {};%
  \end{tikzpicture}%
  }
}

\newcommand{\linestraight}[3]{%
  \raisebox{0.5\height}{
  \begin{tikzpicture}
    \draw[line width=2.5pt, color={rgb,255:red,#1;green,#2;blue,#3}] (0,-.1) -- (0.75,-.1);
  \end{tikzpicture}%
  }
}

\newcommand{\linediamond}[3]{%
  \begin{tikzpicture}[baseline={(0,-.1)}]
    \coordinate (A) at (0.15,0);
    \coordinate (B) at (.25,-.1);
    \coordinate (C) at (.35,0);
    \coordinate (D) at (.25,.1);
    \draw[draw=none, fill={rgb,255:red,#1;green,#2;blue,#3}] (A) -- (B) -- (C) -- (D) -- cycle;
    \draw[line width=0.8pt, color={rgb,255:red,#1;green,#2;blue,#3}] (\startlegend,0) -- (0.6,0);
  \end{tikzpicture}%
  }

\newcommand{\diagonalcross}[3]{%
  \begin{tikzpicture}[baseline={(0,-.1)}]
    \filldraw[fill={rgb,255:red,#1;green,#2;blue,#3}, draw=none, rotate=45] (-0.1,-0.03) rectangle (0.1,0.03);
    \filldraw[fill={rgb,255:red,#1;green,#2;blue,#3}, draw=none, rotate=-45] (-0.1,-0.03) rectangle (0.1,0.03);
    \draw[line width=0.8pt, color={rgb,255:red,#1;green,#2;blue,#3}] (-.45,0) -- (0.45,0);
  \end{tikzpicture}%
}

\usetikzlibrary{shapes.geometric}

\newcommand{\linecircle}[3]{%
  \begin{tikzpicture}[baseline={(0,-.1)}]
    \draw[draw=none, fill={rgb,255:red,#1;green,#2;blue,#3}](0.25, 0) circle (.08);
    \draw[line width=0.8pt, color={rgb,255:red,#1;green,#2;blue,#3}] (\startlegend,0) -- (\enlegend,0);
  \end{tikzpicture}%
  }

\newcommand{\linedotted}[3]{%
\begin{tikzpicture}[baseline={(0,-.1)}]
  \draw[line width=1.5pt, dotted, color={rgb,255:red,#1;green,#2;blue,#3}] (\startlegend,0) -- (\enlegend,0);
\end{tikzpicture}%
}

\newcommand{\linesquare}[3]{%
  \begin{tikzpicture}[baseline={(0,-.1)}]
    \draw[draw=none, fill={rgb,255:red,#1;green,#2;blue,#3}](0.179, -.071) rectangle (.321,.071);
    \draw[line width=0.8pt, color={rgb,255:red,#1;green,#2;blue,#3}] (\startlegend,0) -- (\enlegend,0);
  \end{tikzpicture}%
  }

\newcommand{\linelefttri}[3]{%
  \begin{tikzpicture}[baseline={(0,-.1)}]
    \coordinate (A) at (0.17,0);
    \coordinate (B) at (.3,-.1);
    \coordinate (C) at (.3,.1);
    \draw[draw=none, fill={rgb,255:red,#1;green,#2;blue,#3}] (A) -- (B) -- (C) -- cycle;
    \draw[line width=0.8pt, color={rgb,255:red,#1;green,#2;blue,#3}] (\startlegend,0) -- (\enlegend,0);
  \end{tikzpicture}%
}

\newcommand{\linerighttri}[3]{
  \begin{tikzpicture}[baseline={(0,-.1)}]
    \coordinate (A) at (0.33,0);
    \coordinate (B) at (.2,-.1);
    \coordinate (C) at (.2,.1);
    \draw[draw=none, fill={rgb,255:red,#1;green,#2;blue,#3}] (A) -- (B) -- (C) -- cycle;
    \draw[line width=0.8pt, color={rgb,255:red,#1;green,#2;blue,#3}] (\startlegend,0) -- (\enlegend,0);
  \end{tikzpicture}%
}

\newcommand{\linedowntri}[3]{
  \begin{tikzpicture}[baseline={(0,-.1)}]
    \coordinate (A) at (0.25,-.09);
    \coordinate (B) at (.17,.06);
    \coordinate (C) at (.33,.06);
    \draw[draw=none, fill={rgb,255:red,#1;green,#2;blue,#3}] (A) -- (B) -- (C) -- cycle;
    \draw[line width=0.8pt, color={rgb,255:red,#1;green,#2;blue,#3}] (\startlegend,0) -- (\enlegend,0);
  \end{tikzpicture}%
}

\newcommand{\lineuptri}[3]{
  \begin{tikzpicture}[baseline={(0,-.1)}]
    \coordinate (A) at (0.25,.09);
    \coordinate (B) at (.17,-.06);
    \coordinate (C) at (.33,-.06);
    \draw[draw=none, fill={rgb,255:red,#1;green,#2;blue,#3}] (A) -- (B) -- (C) -- cycle;
    \draw[line width=0.8pt, color={rgb,255:red,#1;green,#2;blue,#3}] (\startlegend,0) -- (\enlegend,0);
  \end{tikzpicture}%
}

\newcommand{\rkfvisualslegend}{\vspace{-1mm} \small \\
\ccircle{253}{186}{187} Noisy measurements \quad
\ccircle{77}{175}{74} Optimal solution \\
\cblock{0}{0}{255} LAH  \quad
\cblock{150}{150}{150} L2WS  \quad
\cblock{255}{189}{8} LM  \quad
}

\newcommand{\legendridge}{\vspace{-1mm} \small \linedowntri{0}{0}{0} vanilla \hspace{1mm}
\diagonalcross{128}{128}{128} conjugate gradient \hspace{1mm}
\lineuptri{255}{127}{0} Nesterov \hspace{1mm}
  \linediamond{77}{175}{74} silver \hspace{1mm}\\
  \linelefttri{166}{86}{40} nearest neighbor \hspace{1mm}
  \linerighttri{152}{78}{163} L2WS \hspace{1mm}
  \linecircle{228}{26}{28} LM 
  \hspace{1mm}
  \linesquare{55}{126}{184} LAH
}

\newcommand{\legendlogistic}{\vspace{-1mm} \small \linedowntri{0}{0}{0} vanilla \hspace{1mm}
\lineuptri{255}{127}{0} Nesterov \hspace{1mm}
  \linediamond{77}{175}{74} silver \hspace{1mm}\\
  \linelefttri{166}{86}{40} nearest neighbor \hspace{1mm}
  \linerighttri{152}{78}{163} L2WS \hspace{1mm}
  \linecircle{228}{26}{28} LM 
  \hspace{1mm}
  \linesquare{55}{126}{184} LAH
  \linedotted{55}{126}{184} 95th quantile LAH
}

\newcommand{\legendlasso}{
\begin{flushleft}
  \begin{minipage}[t]{0.53\linewidth}
    \centering
    \vspace{-1mm} \small \linedowntri{0}{0}{0} Vanilla \hspace{5mm}
    \lineuptri{255}{127}{0} FISTA \hspace{1mm}\\
      \linelefttri{166}{86}{40} nearest neighbor \hspace{1mm}
      \linerighttri{152}{78}{163} L2WS \hspace{1mm}
      \linecircle{228}{26}{28} LM 
      \hspace{1mm}
      \linesquare{55}{126}{184} LAH
      \hspace{1mm}
      \linedotted{55}{126}{184} 95th quantile LAH
      \end{minipage}
  \begin{minipage}[t]{0.45\linewidth}
    \centering
    \legendlassostep
      \end{minipage}
\end{flushleft}
}

\newcommand{\legendlogisticstep}{\vspace{-1mm} \small
  \cblock{30}{30}{30} step-varying \hspace{2mm} 
  \cblock{150}{150}{150} steady-state 
  \linestraight{251}{124}{163} $2 / L$ \hspace{2mm}
}

\newcommand{\legendlassostep}{\vspace{-1mm} \small
  \cblock{30}{30}{30} step-varying \hspace{2mm} \\
  \cblock{150}{150}{150} steady-state \\
  \linestraight{251}{124}{163} $2 / L$ \hspace{2mm}
}

\newcommand{\legendrkf}{\vspace{-1mm} \small \linedowntri{0}{0}{0} vanilla \hspace{1mm}
  \linelefttri{166}{86}{40} nearest neighbor ($N=10000$) \hspace{1mm}
  \linediamond{77}{175}{74} previous solution warm start \\
  \linerighttri{152}{78}{163} L2WS ($N=10000$) \hspace{1mm}
  \linecircle{228}{26}{28} LM ($N=10000$) \hspace{1mm}
  \linesquare{55}{126}{184} LAH ($N=10$) \hspace{1mm}
}

\newcommand{\legend}{\vspace{-1mm} \small \linedowntri{0}{0}{0} vanilla \hspace{1mm}
  \linelefttri{166}{86}{40} nearest neighbor ($N=10000$) \hspace{1mm}
  \linerighttri{152}{78}{163} L2WS ($N=10000$) \hspace{1mm}\\
  \linecircle{228}{26}{28} LM ($N=10000$) \hspace{1mm}
  \linesquare{55}{126}{184} LAH ($N=10$) \hspace{1mm}
  \linedotted{55}{126}{184} 95th quantile LAH
}

\newcommand{\iters}{\small{Mean iterations to reach a given primal and dual residual (Tol.)}}

\newcommand{\itersunconstrained}{\small{Mean iterations to reach a given suboptimality (Tol.)}}

\newcommand{\figsize}{0.8}

\newenvironment{talign*}
 {\let\displaystyle\textstyle\csname align*\endcsname}
 {\endalign}

\newcommand*{\param}{x}

\newcommand{\tableheader}{
  \begin{tabular}{@{}c@{}}Tol. \\\end{tabular}
  &
\begin{tabular}{@{}c@{}}Vanilla\end{tabular}
&
\begin{tabular}{@{}c@{}}Nearest \\ Neighbor\end{tabular}
&
\begin{tabular}{@{}c@{}}L2WS \\ $N=10$\end{tabular}
&
\begin{tabular}{@{}c@{}}L2WS \\ $N=10000$\end{tabular}
&
\begin{tabular}{@{}c@{}}LM \\ $N=10$\end{tabular}
&
\begin{tabular}{@{}c@{}}LM \\ $N=10000$\end{tabular}
&
\begin{tabular}{@{}c@{}}LAH \\ $N=10$\end{tabular}
\\}

\newcommand{\colnames}{\colA & \colB & \colC & \colD &\colE & \colF &\colG & \colH} 

\newcommand{\myCSVReaderIters}[1]{
\tableheader
      \midrule
    \csvreader[
        head to column names,
        late after line=\\
    ]{#1}{ 
        accuracies=\colA,
        cold_start=\colB,
        nearest_neighbor=\colC,
        l2ws=\colD,
        l2ws10000=\colE,
        lm=\colF,
        lm10000=\colG,
        lasco=\colH,
    }{\colnames}
}

\newcommand{\myCSVReaderItersFull}[1]{
\tableheaderfull
      \midrule
    \csvreader[
        head to column names,
        late after line=\\
    ]{#1}{ 
        accuracies=\colA,
        cold_start=\colB,
        nearest_neighbor=\colC,
        prev_sol=\colD,
        l2ws=\colE,
        l2ws10000=\colF,
        lm=\colG,
        lm10000=\colH,
        lasco=\colI,
    }{\colnamesfull}
}

\newcommand{\tableheaderfull}{
  \begin{tabular}{@{}c@{}}Tol. \\\end{tabular}
  &
\begin{tabular}{@{}c@{}}Vanilla\end{tabular}
&
\begin{tabular}{@{}c@{}}Near. Neigh. \\ $N=10000$\end{tabular}
&
\begin{tabular}{@{}c@{}}Prev. Sol.\end{tabular}
&
\begin{tabular}{@{}c@{}}L2WS \\ $N=10$\end{tabular}
&
\begin{tabular}{@{}c@{}}L2WS \\ $N=10000$\end{tabular}
&
\begin{tabular}{@{}c@{}}LM \\ $N=10$\end{tabular}
&
\begin{tabular}{@{}c@{}}LM \\ $N=10000$\end{tabular}
&
\begin{tabular}{@{}c@{}}LAH \\ $N=10$\end{tabular}\\}

\newcommand{\colnamesfull}{\colA & \colB & \colC & \colD  & \colE & \colF & \colG & \colH & \colI}

\newcommand{\myCSVReaderItersRidge}[1]{
\tableheaderRidge
      \midrule
    \csvreader[
        head to column names,
        late after line=\\
    ]{#1}{ 
        accuracies=\colA,
        cold_start=\colB,
        nesterov=\colC,
        conj_grad=\colD,
        silver=\colE,
        nearest_neighbor=\colF,
        l2ws=\colG,
        lm=\colH,
        lasco_one_step=\colI,
        lasco_two_step=\colJ,
        lasco_three_step=\colK,
        lasco=\colL
    }{\colnamesRidge}
}

\newcommand{\tableheaderRidge}{
  \begin{tabular}{@{}c@{}}Tol. \\\end{tabular}
  &
\begin{tabular}{@{}c@{}}Vanilla\end{tabular}
&
\begin{tabular}{@{}c@{}}Nesterov \\ \end{tabular}
&
\begin{tabular}{@{}c@{}}Conjugate \\ gradient\end{tabular}
&
\begin{tabular}{@{}c@{}}Silver \\ \end{tabular}
&
\begin{tabular}{@{}c@{}}Nearest \\ Neighbor\end{tabular}
&
\begin{tabular}{@{}c@{}}L2WS\end{tabular}
&
\begin{tabular}{@{}c@{}}LM \\ \end{tabular}
&
\begin{tabular}{@{}c@{}}LAH \\ $B=1$\end{tabular}
&
\begin{tabular}{@{}c@{}}LAH \\ $B=2$\end{tabular}
&
\begin{tabular}{@{}c@{}}LAH \\ $B=3$\end{tabular}
&
\begin{tabular}{@{}c@{}}LAH \\ $B=10$\end{tabular}
\\
}

\newcommand{\colnamesRidge}{\colA  & \colB & \colC & \colD & \colE &\colF & \colG & \colH & \colI &\colJ &\colK & \colL}

\newcommand{\myCSVReaderItersLogistic}[1]{
\tableheaderLogistic
      \midrule
    \csvreader[
        head to column names,
        late after line=\\
    ]{#1}{ 
        accuracies=\colA,
        cold_start=\colB,
        nesterov=\colC,
        silver=\colD,
        nearest_neighbor=\colE,
        l2ws=\colF,
        l2ws10000=\colG,
        lm=\colH,
        lm10000=\colI,
        lasco_one_step=\colJ,
        lasco=\colK
    }{\colnamesLogistic}
}

\newcommand{\tableheaderLogistic}{
  \begin{tabular}{@{}c@{}}Tol. \\\end{tabular}
  &
\begin{tabular}{@{}c@{}}Vanilla\end{tabular}
&
\begin{tabular}{@{}c@{}}Nesterov \\ \end{tabular}
&
\begin{tabular}{@{}c@{}}Silver \\ \end{tabular}
&
\begin{tabular}{@{}c@{}}Nearest \\ Neighbor\end{tabular}
&
\begin{tabular}{@{}c@{}}L2WS \\ $N=10$\end{tabular}
&
\begin{tabular}{@{}c@{}}L2WS \\ $N=10000$\end{tabular}
&
\begin{tabular}{@{}c@{}}LM \\ $N=10$\end{tabular}
&
\begin{tabular}{@{}c@{}}LM \\ $N=10000$\end{tabular}
&
\begin{tabular}{@{}c@{}}LAH \\ $B=1$\end{tabular}
&
\begin{tabular}{@{}c@{}}LAH \\ $B=10$\end{tabular}
\\
}

\newcommand{\colnamesLogistic}{\colA  & \colB & \colC & \colD & \colE &\colF & \colG & \colH & \colI & \colJ & \colK}

\newcommand{\myCSVReaderItersLasso}[1]{
\tableheaderLasso
      \midrule
    \csvreader[
        head to column names,
        late after line=\\
    ]{#1}{ 
        accuracies=\colA,
        cold_start=\colB,
        nesterov=\colC,
        nearest_neighbor=\colD,
        l2ws=\colE,
        l2ws10000=\colF,
        lm=\colG,
        lm10000=\colH,
        lasco=\colI,
    }{\colnamesLasso}
}

\newcommand{\colnamesLasso}{\colA  & \colB & \colC &\colD & \colE  & \colF & \colG &\colH & \colI} 

\newcommand{\tableheaderLasso}{
  \begin{tabular}{@{}c@{}}Tol. \\\end{tabular}
  &
\begin{tabular}{@{}c@{}}Vanilla\end{tabular}
&
\begin{tabular}{@{}c@{}}FISTA \\ \end{tabular}
&
\begin{tabular}{@{}c@{}}Near.  Neigh. \\ $N=10000$ \end{tabular}
&
\begin{tabular}{@{}c@{}}L2WS \\ $N=10$\end{tabular}
&
\begin{tabular}{@{}c@{}}L2WS \\ $N=10000$\end{tabular}
&
\begin{tabular}{@{}c@{}}LM \\ $N=10$ \end{tabular}
&
\begin{tabular}{@{}c@{}}LM \\ $N=10000$ \end{tabular}
&
\begin{tabular}{@{}c@{}}LAH \\ $N=10$\end{tabular}
\\
}

\newcommand{\myCSVReaderRed}[1]{%
\tableheader
      \midrule
    \csvreader[
        head to column names,
        late after line=\\
    ]{#1}{ 
        accuracies=\colA,
        cold_start_red=\colC,
        nearest_neighbor_red=\colD,
        maml_red=\colO,
        obj_k0_red=\colE,
        obj_k5_red=\colF,
        obj_k15_red=\colG,
        obj_k30_red=\colH,
        obj_k60_red=\colI,
        reg_k0_red=\colJ,
        reg_k5_red=\colK,
        reg_k15_red=\colL,
        reg_k30_red=\colM,
        reg_k60_red=\colN,
    }{\colnames}
}

\newcommand{\myCSVReaderRedAlt}[2]{%
    \tableheaderAlt{#1}
      \midrule
    \ifstrequal{#1}{MAML}{%
    \csvreader[
        head to column names,
        late after line=\\
    ]{#2}{ 
        accuracies=\colA,
        cold_start_red=\colC,
        nearest_neighbor_red=\colD,
        maml_red=\colO,
        obj_k0_red=\colE,
        obj_k1_red=\colF,
        obj_k5_red=\colG,
        obj_k15_red=\colH,
        obj_k60_red=\colI,
        reg_k0_red=\colJ,
        reg_k1_red=\colK,
        reg_k5_red=\colL,
        reg_k15_red=\colM,
        reg_k60_red=\colN,
    }{\colnames}
    }
    {
      \csvreader[
        head to column names,
        late after line=\\
    ]{#2}{ 
        accuracies=\colA,
        cold_start_red=\colC,
        nearest_neighbor_red=\colD,
        prev_sol_red=\colO,
        obj_k0_red=\colE,
        obj_k1_red=\colF,
        obj_k5_red=\colG,
        obj_k15_red=\colH,
        obj_k60_red=\colI,
        reg_k0_red=\colJ,
        reg_k1_red=\colK,
        reg_k5_red=\colL,
        reg_k15_red=\colM,
        reg_k60_red=\colN,
    }{\colnamesprevsol}
    }
}

\newcommand{\useCSVReaderRed}[2]{%
    \ifdefined\currentCSVReaderRed
        \expandafter\csname\currentCSVReaderRed\endcsname{#1}{#2}%
    \else
        \PackageWarning{YourPackage}{No CSV reader command defined, defaulting to \myCSVReaderRed}%
        \myCSVReaderRed{#1}{#2}%
    \fi
}

\makeatletter
\renewcommand{\eqref}[1]{\textup{\tagform@{\ref{#1}}}}
\makeatother

\begin{document}
\maketitle
\begin{abstract}
  We introduce a machine-learning framework to learn the hyperparameter sequence of first-order methods (\eg, the step sizes in gradient descent) to quickly solve parametric convex optimization problems.
  Our computational architecture amounts to running fixed-point iterations where the hyperparameters are the same across all parametric instances and consists of two phases.
  In the first \emph{step-varying} phase the hyperparameters vary across iterations, while in the second \emph{steady-state} phase the hyperparameters are constant across iterations.
  Our learned optimizer is flexible in that it can be evaluated on any number of iterations and is guaranteed to converge to an optimal solution.
  To train, we minimize the mean square error to a ground truth solution.
  In the case of gradient descent, the one-step optimal step size is the solution to a least squares problem, and in the case of unconstrained quadratic minimization, we can compute the two and three-step optimal solutions in closed-form.
  In other cases, we backpropagate through the algorithm steps to minimize the training objective after a given number of steps.
  We show how to learn hyperparameters for several popular algorithms: gradient descent, proximal gradient descent, and two ADMM-based solvers: OSQP and SCS.
  We use a sample convergence bound to obtain generalization guarantees for the performance of our learned algorithm for unseen data, providing both lower and upper bounds.
  We showcase the effectiveness of our method with many examples, including ones from control, signal processing, and machine learning.
  Remarkably, our approach is highly data-efficient in that we only use $10$ problem instances to train the hyperparameters in all of our examples.
\end{abstract}

\ifpreprint \else
\fi

\newcommand{\myparagraph}[1]{%
  \paragraph{#1\ifpreprint.\fi}%
}





\section{Introduction}\label{sec:intro}
This paper considers \emph{parametric fixed-point problems} of the form
\begin{equation}\label{prob:fp}
  \mbox{find} \; z \quad \mbox{ such that } \quad  z = T(z, x),
\end{equation}
where $z \in \reals^\fplen$ is the decision variable, $x \in \reals^d$ is the context or parameter drawn from a distribution $\mathcal{X}$, and $T : \reals^\fplen \times \reals^d \rightarrow \reals^\fplen$ is a mapping.
Problem~\eqref{prob:fp} implicitly defines a (potentially non-unique) solution $z^\star(x)$, which we assume always exists.
We focus exclusively on the setting where problem~\eqref{prob:fp} represents the optimality conditions of a parametric convex optimization problem.
Indeed, nearly all convex optimization problems can be cast as fixed-point problems, which typically represent the optimality conditions~\citep{lscomo}.
In many applications, it is common to repeatedly solve problem~\eqref{prob:fp} but with a varying parameter $x$.
For example, in control and robotics, we repeatedly solve optimization problems as the state changes to update the control inputs~\citep{borrelli_mpc_book}.
This parametric structure also arises in other applications like signal processing where signals are repeatedly recovered from noisy measurements in the same system~\citep{lista} and optimal power flow where the generator outputs and voltage magnitudes are updated in response to varying demand and renewable generation~\citep{ml_opf}.
Solving parametric convex optimization problems are often bottlenecks in the systems they are apart of, and typically, iterative algorithms are needed to solve them.
Because of their cheap per-iteration cost, first-order methods which only use first-order information~\citep{fom_book} and take the form
\begin{equation}\label{eq:fp}
  z^{k+1}(x) = T(z^k(x), x),
\end{equation}
are a popular choice to solve these convex problems.
Under suitable conditions on the operator~$T$, which typically hold for convex optimization, the iterates are known to converge to an optimal solution, \ie, there exists an optimal solution $z^\star(x)$ such that $\lim_{k \rightarrow \infty} \|z^k(x) - z^\star(x)\|_2 = 0$.
Yet, in many applications, the number of fixed-point iterations that can be run is limited due to solve-time constraints.
For example, in many problems in control and robotics, we only have milliseconds before we need to solve the next parametric problem~\citep{borrelli_mpc_book}.
This can lead to suboptimal or infeasible solutions, which may not be acceptable for safety-critical applications.

\myparagraph{Speeding up the convergence of fixed-point iterations}
Despite the widespread use of fixed-point optimization algorithms, they can suffer from slow convergence~\citep{zhang2020globally}.
One approach to improve the convergence speed is to carefully choose the algorithm hyperparameters.
The mainstream approach is to use constant algorithm hyperparameters, and in certain cases, like gradient descent for smooth, convex optimization, the optimal step size has been well-established~\citep{mon_primer,lscomo}.
Recently, interest has grown in \emph{varying} the step sizes across the iterations~\citep{grimmer2024provably,grimmer2023accelerated} for gradient descent.
As an example, the silver step size rule~\citep{altschuler2023acceleration_str_cvx,altschuler_nonstr} judiciously picks non-constant step sizes to improve the convergence rate; however, while it has been numerically proven to lead to faster behavior for constrained parametric quadratic optimization~\cite{perfverifyqp}, its convergence analysis is currently limited to unconstrained optimization problems.

An alternative strategy is to use acceleration methods~\citep{acceleration_survey}, which combine past iterates to generate the next one.
In certain cases, these methods provably improve the worst-case convergence rate, \eg, Nesterov's method~\citep{nesterov} for unconstrained, smooth convex optimization.
Acceleration techniques like Anderson acceleration~\citep{anderson_acceleration} have been developed for more general cases, but devising acceleration techniques that reliably work well remains a challenge.
Most importantly, neither of these two strategies, carefully picking algorithm hyperparameters nor acceleration techniques, takes advantage of the parametric nature of problem~\eqref{prob:fp}.

\myparagraph{Learning to optimize}
In recent years, a paradigm called \emph{learning to optimize} has gained attention for its ability to use machine learning to reduce the solve-time in parametric optimization algorithms.
Techniques include learning initializations~\citep{l2ws}, metrics for proximal algorithms~\citep{metric_learning}, and step sizes~\citep{rlqp}.
Learned optimizers have also been used to solve other types of problems, \eg, inverse problems~\citep{lista,alista}, non-convex optimization~\citep{e2e_survey}, and meta-learning~\citep{maml}.
Yet, several barriers limit the effectiveness of existing learned optimizers.

First, since many learned optimizers replace algorithm steps with learned variants, asymptotic convergence may no longer be guaranteed~\citep{amos_tutorial,l2o}.
While some methods have been able to ensure the convergence of learned optimizers, \eg, by learning warm starts~\citep{l2ws}, safeguarding algorithms steps~\citep{safeguard_convex}, or ensuring that the steps do not deviate too much from a method that is known to converge~\citep{banert2021accelerated}, guaranteeing convergence without subtracting from strong empirical performance persists as a challenge.

Second, learned optimizers typically lack generalization guarantees~\citep{amos_tutorial,l2o}, which can provide numerical guarantees that hold with high probability within a budget of iterations.
Recent work has been used to endow learned optimizers with PAC-Bayes generalization guarantees~\citep{sambharya2024data,Sucker2024LearningtoOptimizeWP}.
Yet, these methods can struggle to provide non-vacuous guarantees, for example when the number of samples is small~\citep{sambharya2024data}.

Third, there is no clear best way to pick the right computational architecture for the learned optimizer.
Developing methods that are robust, versatile, and reliably perform well on a wide-range of parametric families remains an open problem~\citep{amos_tutorial}. 

\myparagraph{Contributions} We propose a framework to learn algorithm hyperparameters (LAH)  of first-order methods in parametric convex optimization problems.
Our detailed contributions are as follows:
\begin{itemize}
  \item We introduce a two-phase architecture where the hyperparameters are shared across all problem instances.
In the first \emph{step-varying} phase, the set of learnable hyperparameters are free to vary across the iterations.
In the second \emph{steady-state} phase, we learn a fixed set of hyperparameters that remain constant from a predetermined number of iterations onwards.
If run for enough iterations, the iterates are guaranteed to converge to an optimal solution due to the inclusion of the steady-state phase in our architecture.
  \item We introduce a progressive training strategy in which we optimize a given number of steps at a time. 
  For gradient descent, we show that the one-step optimal step size is the solution of a least squares problem.
  In the special case of quadratic minimization, we show that the optimal two and three-step optimal step sizes can be computed in closed-form, and moreover, if the parameter $x$ is drawn from a known Gaussian, no training instances are needed, and we can solve the stochastic problem directly.
  In other cases, we use gradient-based methods to optimize the algorithm hyperparameters.
  \item Using sample convergence bounds with validation datasets, we construct generalization guarantees that hold with high probability. 
  We then use these guarantees to construct upper and lower quantile bounds for the performance metric of interest.
  \item 
  We learn the algorithm hyperparameters for several popular methods include gradient descent, proximal gradient descent, and two popular ADMM-based solvers: the operator-splitting solver for quadratic programs (OSQP)~\citep{osqp} and the splitting conic solver (SCS)~\citep{scs_quadratic}.
  \item We showcase the strength of our method across several applications in control, signal processing, and machine learning. We compare against existing methods for learned optimizers, namely a method that learns the initialization~\citep{l2ws} and one that learns the metric from data~\citep{metric_learning}.
  Notably, our approach is highly data-efficient in that we only use $10$ training samples in each example to train our weights.
\end{itemize}

\myparagraph{Notation}
We denote the set of vectors with real entries, non-negative real entries, and positive entries, each of size $n$, with $\reals^n$, $\reals^n_+$, and $\reals^n_{++}$ respectively.
We denote the set of symmetric matrices, positive semidefinite matrices, and positive definite matrices, each of size $n \times n$, with $\symm^n$, $\symm^n_+$, and $\symm^n_{++}$ respectively.
We use $\mathbf{E}$ to denote expectation.
The trace of a matrix $A$ is given by $\Tr(A)$.
Given two vectors $u \in \reals^n$ and $v \in \reals^n$, $u \odot v$ denotes their element-wise product.
We write the vector of all ones with length $d$  as $\mathbf{1}_d$ and the identiy matrix of shape $d \times d$ as $I_d$.
For a boolean conditions $c$, we let $\mathbf{1}(c) = 1$ if $c$ is true, and $0$ otherwise.
For a matrix $A \in \reals^{m \times n}$ we denote its largest singular value as $\|A\|_2 =\max_{\|u\|_2 = 1} \|A u\|_2$.
We denote the proximal operator $\prox_h : \reals^n \rightarrow \reals^n$ as~\citep{prox_algos} $\prox_h(v) = \argmin_u h(u) + (1 / 2) \|u - v\|_2^2$.
For a vector $v \in \reals^N$, the geometric mean is computed as $(\Pi_{i=1}^N v_i)^{1 / N}$.

\section{Related work}
We now describe various areas in which learning for optimization has recently gained wide attention.
\myparagraph{Learning algorithm steps for convex optimization}
Several works aim to learn the key components of the algorithm steps for convex optimization.
In the convex QP case, reinforcement learning can accelerate the OSQP solver by dynamically predicting the step size~\citep{rlqp}.
More generically in convex optimization, learning via stochastic gradient descent can accelerate fixed-point iterations via recurrent neural networks~\citep{neural_fp_accel_amos} or by tuning the metric used in proximal algorithms~\citep{metric_learning}, with, in some cases, provable convergence guarantees~\citep{banert2021accelerated}.
Our work shares the same aim as these works, but is different in method.
By only learning a shared algorithm hyperparameter sequence, we significantly reduce the number of weights to learn and, therefore, the number of samples required, making it highly data-efficient.
In all of our numerical examples, we only use ten training instances.
In contrast, in all of the works mentioned above, thousands of training instances are needed.
It is important to note that our method is supervised as it requires access to solutions of the training problem instances, as opposed to some other techniques which rely on unsupervised formulations~\citep{neural_fp_accel_amos,rlqp}.
An additional advantage of our approach, specific to the case of solvers requiring linear system solutions (\eg, OSQP and SCS), is that we do not need to refactor any matrices to solve new problems.
Instead, other works~\citep{metric_learning,rlqp} require at least one matrix factorization.

\myparagraph{Learning algorithm steps beyond convex optimization}
The idea of learning algorithm steps has also been applied to inverse problems where the goal is to recover a true signal rather than to minimize an objective.
This strategy has been applied to a variety of domains such as sparse coding~\citep{lista,alista}, image restoration~\citep{Diamond2017UnrolledOW,plug_and_play_ryu}, and wireless communication systems~\citep{deep_unfolding_wireless}.
Learning algorithm steps has also been applied to solve non-convex optimization problems, \eg, by designing algorithms for combinatorial problems~\citep{balcan2020data}, and learning branching heuristics for integer programming~\citep{learn_2_branch}.
Our work differs in scope from these methods since we focus on algorithms to solve convex optimization problems.

\myparagraph{Learning initializations}
For convex optimization, instance-specific initializations have been learned for fixed-point problems~\citep{l2ws} and quadratic programs~\citep{l2ws_l4dc} that are tailored for the downstream algorithm.
Other works learn initializations from data, typically in a decoupled fashion, for example in optimal power flow~\citep{warm_start_power_flow} and trajectory optimization~\citep{constraint_informed_traj_ws,li2024efficient}.
In model predictive control~\citep{borrelli_mpc_book} where similar problems need to be repeatedly solved, machine-learning methods have been explored to warm-start an active set method~\citet{mpc_primal_active_sets}, and project the output of a neural network prediction onto the feasible region~\citep{deep_learning_mpc_karg,explicit_mpc_neural_net}.
The approach of learning initializations also been used in meta-learning where a shared initial weights is learned~\citep{maml}, and in non-convex optimization where multiple candidate initial point are learned~\citep{sharony2024learning}.
We find that our method is more data-efficient than this approach to learn the initialization (we compare numerically against~\citep{l2ws}) and that it can handle a wider distribution of problems.
Moreover, in general, the learned initializations may struggle to effect the tail convergence rate.

\myparagraph{Algorithm design}
Choosing the algorithm hyperparameters (\eg, the step size in gradient descent) has been an intensely-studied aspect of designing first-order methods.
The mainstream approach is to use a constant hyperparameter set.
In the case of using gradient descent to solve convex, unconstrained optimization problems, the optimal step size (in terms of worst-case guarantee) is a function of the smoothness parameter of the objective~\citep{mon_primer}.
In the case of a strongly convex objective, the optimal step size is a function of the smoothness and strongly convex parameters~\citep{mon_primer}.
For the more general constrained case, a method based on semidefinite programming has been used to provide the optimal metric for ADMM in the case where the objective can be split into the sum of two parts, one of them being strongly convex and smooth~\citep{giselsson_lin_conv}.

It is natural to wonder if a varying step size schedule can improve convergence.
For the case of unconstrained quadratic minimization, Young's Chebyshev step sizes were shown to improve the convergence rate of gradient descent in 1953~\citep{Young1953}.
In the past few years, several works have shown that varying the step sizes in the more general case of smooth, convex optimization can provably converge faster~\citep{grimmer2023accelerated,grimmer2024provably}.
For example, the silver step size rule~\citep{altschuler2023acceleration_str_cvx,altschuler_nonstr} improves the convergence rate compared with a constant step size schedule.
In the more general setting of applying ADMM to solve constrained, convex optimization problems (\eg, quadratic and conic programs), it has been proven to be a more challenging task to choose the best algorithm hyperparameters.
The performance estimation problem framework~\citep{pep} has been used to construct a non-convex optimization problem to design optimal first-order methods~\citep{pesto_bnb}.
However, none of these methods leverage the parametric structure of our interest.
In this paper, we take advantage of this structure and learn algorithm hyperparameters that are allowed to vary across iterations.
With this approach, we can drastically improve the convergence speed of first-order methods.




\myparagraph{Generalization guarantees for learned optimizers}
Despite the strong empirical performance of learned optimizers, these methods typically lack generalization guarantees to unseen data~\citep{amos_tutorial,l2o}.
To address this limitation, several works have developed methods to provide such guarantees for learned optimizers, for example by using the PAC-Bayes framework~\citep{sambharya2024data,Sucker2024LearningtoOptimizeWP}.
Yet, these methods can struggle when the number of samples is low~\citep{sambharya2024data}.
On the other hand, generalization guarantees that are derived from a validation set are known to be stronger numerically~\citep{kawaguchi2017generalization}.
Motivated by this, we construct numerically strong guarantees using a validation dataset, with only $1000$ instances.

\section{Learning algorithm hyperparameters framework}\label{sec:framework}
In this section we present our computational architecture based on a generic fixed-point algorithm with weights that correspond to the hyperparameters of the algorithm.
We first explain how to run LAH given weights $\theta$, then how to evaluate LAH, and finally how train LAH to minimize a performance loss.

\myparagraph{Running LAH}
Running LAH amounts to running fixed-point iterations with different hyperparameters across the iterations.
The weights are $\theta = (\theta^0, \dots, \theta^H)$ and each $\theta^k \in \reals^a$ corresponds to $a$ hyperparameters of the algorithm used at the $k$-th iteration.
Thus $\theta$ is a vector in $\reals^{(H+1)a}$.
Our computational architecture consists of two phases.
\begin{itemize}
  \item {\bf Step-varying phase}. 
  The first phase is the step-varying phase, where the $k$-th iterate uses the hyperparameters $\theta^k$.
  This phase holds until $H$ steps have been executed.
  The iterates in this phase can be written as
  \begin{equation*}
    z^{k+1}_\theta(x) = T_{\theta^k}(z^k_\theta(x), \param),\quad k=0,1,\dots, H-1.
  \end{equation*}
  \item {\bf Steady-state phase}. The second phase is the steady-state phase which uses the steady-state algorithm hyperparameters $\theta^H$ after $K$ steps:
  \begin{equation*}
    z^{k+1}_\theta(x) = T_{\theta^H}(z^k_\theta(x), \param),\quad k \geq H.
  \end{equation*}
\end{itemize}
The step-varying phase is the core part of our method.
We find, that by allowing the hyperparameters to vary across iterations (as opposed to following the mainstream approach of keeping them fixed), we can improve the performance by significant margins.

There are several reasons why we choose to include a steady-state phase in our architecture where hyperparameters stop  varying. 
First, it allows us to guarantee convergence by relying on existing theory for first-order methods with fixed hyperparameters.
As long as the steady-state hyperparameters $\theta^H$ are appropriately chosen (\eg, they are positive or lie within a range), LAH will converge to an optimal solution regardless of the hyperparameters in the step-varying phase.
Second, by having a fixed number of hyperparameters to learn, we can precompute the most expensive operations depending on them. 
For example, solvers such as OSQP~\citep{osqp} and SCS~\citep{scs_quadratic}, require a linear system to be solved at each iteration, and precomputing a fixed, finite number of factorizations can greatly speedup the online execution. 
We provide more details about factorization caching in \Sec~\ref{sec:fom}.
We remark that the steady-state phase could be replaced with any algorithm that is guaranteed to converge (\eg, Nesterov's acceleration for gradient descent).

\myparagraph{Evaluating LAH}
To evaluate the performance of LAH for a problem with parameter $x$ at a given point $z$, we use a performance metric $\phi(z,x)$.
In this paper, we use two different performance metrics depending on whether the problem is constrained or not.
\begin{itemize}
  \item {\bf Unconstrained case.} Consider the case where the objective is $f(z,x)$ for a function $f : \reals^\fplen \times \reals^d \rightarrow \reals$ that is convex in $z$.
  Here, our performance metric is the suboptimality $\phi(z,x) = f(z, x) - f(z^\star(x), x)$.
  \item {\bf Constrained case.} For constrained optimization problems, it is common to use the primal and dual residuals~\citep{osqp,scs_quadratic} to measure the performance of algorithms.
  Our performance metric $\phi(z,x)$ for constrained problems is the maximum of the primal and dual residual.
\end{itemize}

\myparagraph{The LAH training problem}
Our training performance loss is the mean square error to a ground truth solution after a chosen number of steps $K$.
The problem formulation is
\begin{equation}\label{prob:simplified_l2o}
  \begin{array}{ll}
\mbox{minimize} & \mathbf{E}_{x \sim \mathcal{X}} \|z_\theta^{K}(x) - z^\star(x)\|_2^2\\
  \mbox{subject to} &z^{k+1}_\theta(x) = T_{\theta^k}(z^k_\theta(x), \param),\quad k=0,1,\dots, H-1 \\
  &z^{k+1}_\theta(x) = T_{\theta^H}(z^k_\theta(x), \param),\quad k=H,H+1,\dots, K \\
  &z^0_\theta(x) = 0.
\end{array}
\end{equation}
The first $H$ steps correspond to the step-varying phase (in the first line of constraints), and steps $H$ through $K$ correspond to the steady-state phase where the hyperparameters $\theta^H$ are the same across the iterations (in the second line of constraints).
We always set the initial point to the zero vector, but in principle, this could be learned as well~\citet{l2ws}.
Since in general we do not have access to the distribution $\mathcal{X}$, we minimize problem~\eqref{prob:simplified_l2o} over $N$ independent and identically distributed (i.i.d.) samples $\{x_i\}_{i=1}^N$.
We remark that the loss function in problem~\eqref{prob:simplified_l2o} requires access to (potentially non-unique) ground-truth solutions $\{z^\star(x_i)\}_{i=1}^N$, and falls under the category of regression-based losses~\citet{amos_tutorial}.



\begin{figure}[!h]
  \centering
    \includegraphics[width=1.0\linewidth]{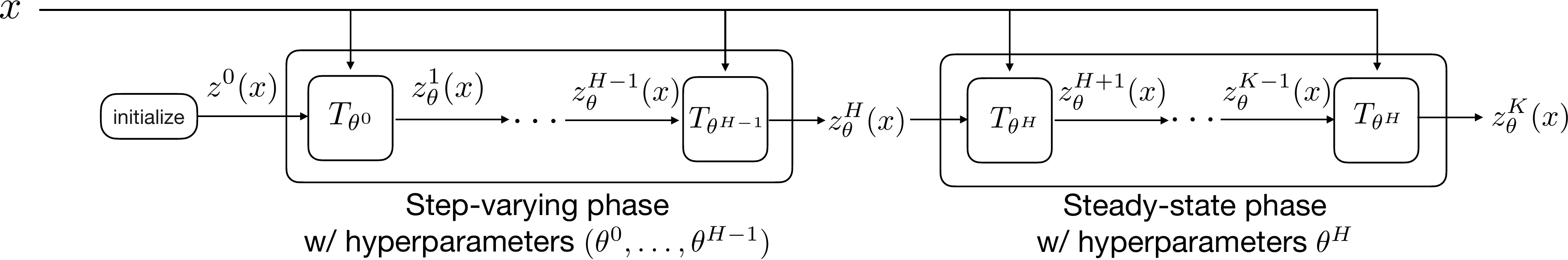}
    \\
    \caption{LAH diagram.
    Running LAH amounts to running fixed-point iterations with the learned hyperparameters $\theta$ across the iterations and consists of two phases.
    First, in the time-varying phase, we run $H$ fixed-point iterations each with a varying hyperparameter set $\theta^k$ in the $k$-th iteration and initialized with $z^0(x) = 0$.
    Second, after the initial $K$ steps, in the steady-state phase, we run $K - H$ fixed-point iterations each using the same hyperparameter set $\theta^H$.
    At evaluation time, we are free to run any number of iterations: not just the number of steps trained on.
    }
    \label{fig:lasco_diagram}
\end{figure}

\myparagraph{Progressive training}
Instead of directly trying to solve problem~\eqref{prob:simplified_l2o} with end-to-end learning techniques which can be difficult to optimize, we apply a progressive training approach, a common strategy in the learning to optimize literature~\citep{l2o,cai2021learned}.
This method breaks problem~\eqref{prob:simplified_l2o} into more manageable pieces by sequentially optimizing a given number of algorithm hyperparameters at a time. 
A key piece of our progressive training procedure is the $B$-step lookahead problem.
Given starting points  $\{\tilde{z}^k(x_i)\}_{i=1}^{N}$, the $B$-step lookahead problem is formulated as
\begin{equation}\label{prob:B_step}
  \begin{array}{ll}
\mbox{minimize} & \sum_{i=1}^N \|z^{B}_\theta(x_i) - z^\star(x)\|_2^2\\
  \mbox{subject to} &z^{l+1}_\theta(x_i) = T_{\theta^l}(z^l_\theta(x_i), x_i),\quad l=k,\dots, k+B-1,\quad i=1, \dots, N \\
  &z^k_\theta(x_i) = \tilde{z}^k(x_i),
\end{array}
\end{equation}
where the decision variables are $\theta^k,\dots,\theta^{k+B-1}$.
The progressive training procedure alternates between the following two steps that i) solve the $B$-step lookahead problem to find the next $B$ hyperparameters, and then ii) subsequently run $B$ steps of the algorithm for each training problem to get the new initial points:
\begin{itemize}
  \item Set $\tilde{z}^k(x_i) = z^k_\theta(x_i)$ and solve problem~\eqref{prob:B_step} to get $(\theta^k, \dots, \theta^{k+B-1})$.
  \item Run $B$ steps: $z^{k+B}_\theta(x_i) = T_{\theta^{k+B-1}}(\cdots T_{\theta^{k+1}} (T_{\theta^k}(z^k_\theta(x_i),x_i),x_i),x_i)$ for $i=1, \dots, N$.
\end{itemize}
In \Sec~\ref{sec:solve_train}, we show how to solve the $B$-step lookahead problems.



\section{First-order methods and their hyperparameters}\label{sec:fom}
In this section, we enumerate several first-order methods and the hyperparameters that are learned in each one in Table~\ref{table:fp_algorithms}.
Each algorithm takes the form of fixed-point iterations from~\eqref{eq:fp}.
In each example, we show how the steady-state hyperparameter $\theta^H$ can be picked to guarantee convergence of our method.
We provide more details in Appendix~\ref{sec:fom_details}.

\begin{table}[!h]
  \centering
  \caption{Several popular first-order methods and their hyperparameters.}
  \label{table:fp_algorithms}
  \adjustbox{max width=\textwidth}{
  \begin{threeparttable}
  \begin{tabular}{@{}lllll@{}}
    \toprule[\heavyrulewidth]
    Algorithm & Problem &  Iterates $z^{k+1} = T(z^k,x)$ & Hyperparameters & \makecell[l]{Steady-state \\ constraints }\\
    \midrule[\heavyrulewidth]
    \makecell[l]{Gradient \\ descent } & $\begin{array}{@{}ll}\min &f(z,x)\end{array}$ & $z^{k+1} = z^k - \theta^k \nabla f(z^k,x)$  & $\theta^k$ & $0 < \theta^H < 2 / L$\\
    \midrule
    \makecell[l]{Proximal \\ gradient \\ descent }
    & $\begin{array}{@{}ll}\min & f(z,x)+ g(z,x)\end{array}$ & $z^{k+1} = \prox_{\theta^k g} (z^k - \theta^k \nabla f(z^k,x), x)$ & $\theta^k$ & $0 < \theta^H < 2 / L$\\
    \midrule
    \makecell[l]{OSQP \\ \citep{osqp}}\textcolor{black}{} & $\begin{array}{@{}ll}
    \min & (1/2) w^T Pw + c^T w\\
    \mbox{s.t.}& l \leq Aw \leq u \quad \text{dual\ } (y) \\
    &\\&\\
    \text{with}&x = (c, l, u)\\
    \end{array}$
    &
    $\begin{aligned}
      & (w^k, \xi^k)= z^k\\
      & v^{k+1} = \Pi_{[l,u]}(\xi^k)\\
      & \text{solve } Q w^{k+1} = \sigma^k w^k - c + \rho^k \odot (A^T (2 v^{k+1} - \xi^k))\\
      & \xi^{k+1} = \rho^k \odot (A w^{k+1} + \xi^{k} - 2 v^{k+1}) + \xi^k - v^{k+1} \\
      & z^{k+1} = (w^{k+1}, \xi^{k+1})\\
      \\
      & \text{with}\quad Q = P + \sigma^k I + \diag (\rho^k)  A^T A\\
      &\rho^k = (\rho^k_{\rm eq} \mathbf{1}_{m_{\rm eq}}, \rho^k_{\rm ineq} \mathbf{1}_{m_{\rm ineq}})\\
    \end{aligned}$ & $\theta^k = (\sigma^k, \rho^k_{\rm eq}, \rho^k_{\rm ineq}, \alpha^k)$ & $\begin{aligned}
      \sigma^H > 0 \\
      \rho^H_{\rm eq} > 0 \\
      \rho^H_{\rm ineq} > 0 \\
      1 < \alpha^H < 2
    \end{aligned}$\\
    \midrule[\lightrulewidth]
    \makecell[l]{SCS \\ \citep{scs_quadratic}}\textcolor{black}{} & $\begin{array}{@{}ll}
    \min & (1/2) w^T Pw + c^T w\\
    \mbox{s.t.}& Aw + s = b \quad \text{dual\ } (y) \\
    & s \in \mathcal{K} \\
    &\\
      \text{with}&x = (c, b)
    \end{array}$ & $\begin{aligned}
      & (\mu^k, \eta^k)= z^k\\
      & \tilde{u}^{k+1} :
      \begin{cases}
        \text{solve } (R^k + M) p^k = \mu^k \\
        \tilde{\tau}^{k+1} = {\bf root_+}(\mu^k, \eta^k, p^k, x, r^k_\tau) \\
        \tilde{w}^{k+1} = p^k - \xi \tilde{\tau}^{k+1}
      \end{cases}\\
      & u^{k+1} :
      \begin{cases}
        w^{k+1} = \Pi_{ \reals^q \times \mathcal{K}}(2 \tilde{w}^{k+1} - \mu^k) \\
        \tau^{k+1} = \Pi_{\reals_+}(2 \tilde{\tau}^{k+1} - \eta^k)
      \end{cases}\\
      & z^{k+1} = (\mu^k + \alpha^k (w^{k+1} - \tilde{w}^{k+1}), \eta^k + \alpha^k (\tau^{k+1} - \tilde{\tau}^{k+1}))\\
      \\& \text{with}\quad 
      R^k = \diag(r_w^k \mathbf{1}_q, r_{y_{\rm z}}^k \mathbf{1}_{m_{\rm z}}, r_{y_{\rm nz}}^k  \mathbf{1}_{m_{\rm nz}}) \\
      &     M = \begin{bmatrix}
        I_q + P & A^T  \\
        A & -I_m \\
      \end{bmatrix}, \xi = (I + M)^{-1}x\\
    \end{aligned}$ & $\theta^k=(r_w^k, r_{y_{\rm z}}^k, r_{y_{\rm nz}}^k, r_\tau^k, \alpha^k)$ & $\begin{aligned}
      r_w^H > 0 \\
      r_{y_{\rm z}}^H > 0 \\
      r_{y_{\rm nz}}^H > 0\\
      r_\tau^H > 0\\
      1 < \alpha^H < 2
    \end{aligned}$
    \\
    \bottomrule[\heavyrulewidth]
  \end{tabular}
  \begin{tablenotes}
  \item  We denote~\prox{} as the proximal operator~\citep{prox_algos}\textcolor{black}{}  (see the notation paragraph in \Sec~\ref{sec:intro}\textcolor{black}{} for a formal definition).
  The notational dependence of the iterates on the parameter $x$ has been removed for simplicity.
  \end{tablenotes}
  \end{threeparttable}
  }
\end{table}

\myparagraph{Caching matrix factorizations}
Running vanilla OSQP or SCS requires factoring a matrix once, and then solving a linear system in each iteration using the same factorization.
Factorizing a dense $n \times n$ matrix has complexity $\mathcal{O}(n^3)$, but solving the linear system provided the factorization has a much cheaper complexity of $\mathcal{O}(n^2)$~\citep[\Sec~11]{vmls}.
In this paper, we are interested in parametric optimization, but only consider the case where the matrices $P$ and $A$ in the quadratic and conic programs in Table~\ref{table:fp_algorithms} are the same for all instances.
Hence, the matrix in the linear system is fixed for all problem instances for vanilla OSQP and SCS, and no online factorizations are required as the same factorization can be used for all problems.

To run LAH, the learned penalty parameters change across iterations, and these hyperparameters change the matrix in the linear system.
Yet, still, LAH does not require any factorizations online.
Once training is complete and we have a set of weights $\theta$, we can do all of the matrix factorizations offline.
We need $K+1$ factorizations ($K$ for the step-varying phase and one for the steady-state phase) to run LAH, but these factorizations can be re-used for all problem instances. 
In contrast, if the hyperparameters are predicted from the problem parameter $x$, as in~\citet{metric_learning}, then there is no way to avoid an online matrix factorization.
If the hyperparameters not only depend on the parameter but also the current iterate as in~\citet{rlqp}, then a matrix is required to be factored in each iteration.

\section{Solving the training problem}\label{sec:solve_train}
In this section, we show how to solve the progressive training subproblems. 
For gradient descent, we show that the one-step lookahead problem can be formulated as a least squares problem in \Sec~\ref{subsec:theory:gd}.
Then, we show that the multi-step lookahead problem in the case of quadratic minimization can be solved to global optimality in \Sec~\ref{subsec:multi_step}.
Next, for quadratic minimization, we show that if the parameter is drawn from a Gaussian distribution, we can directly solve the test problem in \Sec~\ref{subsec:stoch}.
Outside of these cases, we use gradient-based methods to solve the training problem as shown in \Sec~\ref{subsec:grad_method}.
Finally in \Sec~\ref{subsec:safeguard}, we describe a safeguarding mechanism to ensure strong performance of the learned hyperparameter sequence.


\subsection{One-step lookahead analysis for gradient descent}\label{subsec:theory:gd}
In this subsection, we analyze gradient descent.
Consider the progressive training regime where we optimize one step at a time.
For gradient descent, at the $k$-th step, we have access to the iterates $\{z_\theta^k(x_i)\}_{i=1}^N$, and the gradients $\{\nabla f(z^k_\theta(x_i),x_i)\}_{i=1}^N$.
To find the best step size in the $k$-th iteration, we solve the optimization problem
\begin{equation}
  \begin{array}{ll}
  \label{prob:gd_train_single}
  \mbox{minimize} & (1 / N) \sum_{i=1}
  ^N \| z^k_\theta(x_i) - \theta^k \nabla f(z^k_\theta(x_i), x_i) - z^\star(x_i)\|_2^2,
  \end{array}
\end{equation}
where $\theta^k$ is the scalar decision variable and the initial points $z_\theta^k(x_i)$ are known.
Assuming that there exists at least one $z_\theta^k(x_i)$ such that $\nabla f(z_\theta^k(x_i),x_i) \neq 0$ (otherwise all training problems have already been solved), this problem is a least squares problem, whose solution is
\begin{equation}\label{eq:ls_sol}
  \theta^k = \frac{\sum_{i=1}^N \nabla f(z^k_\theta(x_i),x_i)^T (z^k_\theta(x_i) - z^\star(x_i))}{\sum_{i=1}^N \|\nabla f(z^k_\theta(x_i),x_i)\|_2^2}.
\end{equation}
This follows from the first-order optimality conditions of differentiable, convex functions~\citep[\Sec~3.1.3]{cvxbook}.
With the following theorem, we prove that these one-step optimal step sizes as calculated in \Eqn~\eqref{eq:ls_sol} are always non-negative, ensuring that the updates move in the descent direction (as opposed to moving in an ascent direction).
\begin{theorem}\label{proof:pos_step}
  The one-step optimal step size $\theta^k$ for problem~\eqref{prob:gd_train_single} is  non-negative for any possible values of $\{z^k_\theta(x_i)\}_{i=1}^N$ as long as there exists some $z^k_\theta(x_i)$ such that $\nabla f(z^k_\theta(x_i),x_i) \neq 0$.
\end{theorem}
See Appendix \Sec~\ref{sec:pos_step_proof} for the proof.

\subsection{Multi-step lookahead analysis for quadratic minimization}\label{subsec:multi_step}
In subsection~\ref{subsec:theory:gd}, we showed that the one-step lookahead problem admits a closed-form solution in the case of gradient descent.
While progressively training with $B=1$ (\ie, by solving a sequence of least squares problems) offers a tractable way to learn the step size schedule, it does not take advantage of the fact that multiple steps can be optimized simultaneously, which has the potential to improve the performance of our method.
In this section, we show that the two-step and three-step lookahead problems can be solved in closed-form in the case of unconstrained quadratic minimization given by
\begin{equation*}
  \begin{array}{ll}
  \label{prob:gd_unconstrained_quad}
  \mbox{minimize} & (1 / 2) z^T P z + x^T z,
  \end{array}
\end{equation*}
where $z \in \reals^n$ is the decision variable, $P \in \symm_{++}^n$ is fixed across problem instances, and $x \in \reals^n$ is the parameter.
The $B$-step lookahead problem can be formulated as 
\begin{equation}
  \begin{array}{ll}
  \label{prob:multi_step}
  \mbox{minimize} & r(\theta^k, \dots, \theta^{k+B-1}) = \sum_{i=1}^N \|(I - \theta^{k+B-1} P)\cdots(I - \theta^{k} P)(\tilde{z}^k(x_i) - z^\star(x_i))\|_2^2.\!\!\!\!\!\!\!\!\!
  \end{array}
\end{equation}
As mentioned in \Sec~\ref{sec:framework}, the starting points are chosen as $\tilde{z}^k(x_i) = z^k_\theta(x_i)$.
We use the eigendecomposition $P = Q \diag(\lambda) Q^T$ where the eigenvalues $\lambda_1,\dots,\lambda_n$ are all positive.
By using the orthonormal property of matrix $Q$, problem~\eqref{prob:multi_step} can be written as
\begin{equation}
  \begin{array}{ll}
  \label{prob:multi_step_eigen}
  \mbox{minimize} & r(\theta^k, \dots, \theta^{k+B-1})  = \sum_{j=1}^n (1 - \theta^k \lambda_j)^2\cdots (1 - \theta^{k+B-1} \lambda_j)^2 \bar{z}_j^k,
  \end{array}
\end{equation}
where $\bar{z}_j^k = \sum_{i=1}^N (Q^T (z_\theta^k(x_i) - z^\star(x_i)))_j^2 \in \reals_+$ for $j=1,\dots,n$.
Note that in this particular case, the objective value is unaffected by re-ordering the step sizes~\citep{Young1953}. 
In the following theorem, we prove that under some mild conditions, a necessary condition for $(\theta^k,\dots,\theta^{k+B-1})$ to be a local minimizer of problem~\eqref{prob:multi_step} is that the step sizes are all distinct.
\begin{theorem}\label{thm:diff_step_sizes}
  Assume that there exists at least $B$ values of $\bar{z}^k \in \reals^n$ that are non-zero whose corresponding eigenvalues (by index) are all distinct.
  The vector $(\theta^k,\dots,\theta^{k+B-1})$ cannot be a local minimizer of problem~\eqref{prob:multi_step} if there exists $\theta^j = \theta^l$ for $j \neq l$, $j,l \in \{k, \dots, k+B-1\}$.
\end{theorem}
See Appendix~\ref{proof:diff_step_sizes} for the proof.
This observation that there is a disadvantage to a constant step size schedule is in line with Young's varying step sizes for quadratic minimization~\citep{Young1953}.
This case of quadratic minimization is special for another reason; it allows us to obtain a \emph{convergence rate} over \emph{any} parameter $x$, including ones not drawn from the distribution $\mathcal{X}$.
\begin{theorem}[Convergence rate for quadratic minimization~\citep{Young1953}]\label{thm:conv_rate}
  For any vector $x \in \reals^n$ the following convergence rate holds for any sequence $(\theta^0, \dots, \theta^{K-1})$:
  \begin{equation*}
    \|z_\theta^K(x) - z^\star(x)\|_2 \leq \max_j |\Pi_{k=0}^{K-1}(1 - \theta^k \lambda_j)| \|z_\theta^0(x) - z^\star(x)\|_2.
  \end{equation*}
\end{theorem}
Even though we optimize for the case where $x$ is drawn from a distribution $\mathcal{X}$ and fix the initial point $z^0_\theta(x)$ to be zero, we can still obtain a convergence rate on the distance to optimality for any parameter $x$ and any initial point.

\subsubsection{Two-step lookahead analysis for quadratic minimization}\label{sec:two_step}
We first investigate the two-step lookahead problem for quadratic minimization.
We wish to find the step sizes $\theta^{k}$ and $\theta^{k+1}$ to minimize the average distance to optimality after $k+2$ iterations.
The two-step lookahead problem is
\begin{equation}
  \begin{array}{ll}
  \label{prob:two_step}
  \mbox{minimize} & \sum_{i=1}^N \|(I - \theta^{k+1} P) (I - \theta^{k} P)(z_\theta^k(x_i) - z^\star(x_i))\|_2^2,
  \end{array}
\end{equation}
where $\theta^k$ and $\theta^{k+1}$ are the scalar decision variables.
The objective function of problem~\eqref{prob:two_step} is a non-convex quartic polynomial.
Yet, with the following theorem, we show that this problem can be solved in closed-form to global optimality.

\begin{theorem}\label{thm:two_step}
  Assume that there exists at least $2$ values of $\bar{z}^k \in \reals_+^n$ that are non-zero whose corresponding eigenvalues (by index) are distinct.
  Given the quantities
  \begin{equation*}
    a = \sum_{j=1}^n \lambda_j \bar{z}_j^k, \quad b = \sum_{j=1}^n \lambda_j^2 \bar{z}_j^k, \quad c = \sum_{j=1}^n \lambda_j^3 \bar{z}_j^k, \quad d = \sum_{j=1}^n \lambda_j^4 \bar{z}_j^k,
  \end{equation*}
the two-step optimal step sizes are non-negative and given by the expression
  \begin{equation*}
    \frac{ad - bc \pm \sqrt{(bc-ad)^2 - 4 (ac - b^2) (bd - c^2)}}{2 (bd -c^2)},
  \end{equation*}
  in either order.
\end{theorem}
See Appendix \Sec~\ref{proof:two_step} for the proof.
By virtue of \Thm~\ref{thm:diff_step_sizes}, we know that the optimal step sizes (which we argue must exist in the proof) cannot be equal.
This fact allows us to simplify the analysis for unequal step sizes, and in the end show that the optimal step sizes are solutions to a quadratic equation that always has two real roots.
We reiterate that the optimal step sizes are not equal and the order does not matter.
The fact that the step sizes are provably positive ensures that a descent direction is taken in each iteration for each problem.

\subsubsection{Three-step lookahead analysis}
The three-step lookahead problem is
\begin{equation}
  \begin{array}{ll}
  \label{prob:three_step}
  \mbox{minimize} & \sum_{i=1}^N \|(I - \theta^{k+2} P) (I - \theta^{k+1} P)(I - \theta^k P)(z_\theta^k(x_i) - z^\star(x_i))\|_2^2,
  \end{array}
\end{equation}
where $\theta^k$, $\theta^{k+1}$, and $\theta^{k+2}$ are the scalar decision variables.
In the following theorem, we show that we can efficiently solve this problem by finding the roots of a cubic polynomial.
\begin{theorem}\label{thm:three_step}
  Assume that there exists at least $3$ values of $\bar{z}^k \in \reals^n$ that are non-zero whose corresponding eigenvalues by index are distinct.
  Then, the three-step optimal step sizes can be found by finding the roots of a cubic equation.
  Any permutation of those roots is a solution, and the roots are all real-valued.
\end{theorem}
See Appendix \Sec~\ref{proof:three_step} for the proof and for the coefficients of the cubic equation.
A key part of our proof relies on the fact that given the first step $\theta^{k}$, the next two steps can be found using \Thm~\ref{thm:two_step}: \ie, the optimal $\theta^{k+2}$ and $\theta^{k+1}$ are functions of $\theta^k$.


\subsection{Stochastic lookahead analysis for quadratic minimization}\label{subsec:stoch}
When the parameters $x$ are drawn from a known Gaussian distribution, the analysis in this section can be used to directly optimize the multi-step lookahead \emph{test problem}. 
The test problem is a stochastic optimization problem~\citep[Chapter~1]{shapiro2009lectures} as it involves optimizing the expected mean square over the distribution of parameters $\mathcal{X}$.
We let the parameter distribution be $\mathcal{X} = \mathcal{N}(\mu, \Sigma)$, and the $B$-step lookahead problem at the $k$-th iteration boils down to
\begin{equation}
  \begin{array}{ll}\label{prob:multi_step_stochastic}
  \mbox{minimize} & r(\theta^k,\dots,\theta^{k+B-1}) = \underset{x\sim\mathcal{X}}{\mathbf{E}} \|(I - \theta^{k+B-1} P)\cdots(I - \theta^{k} P)(z_\theta^k(x) - z^\star(x))\|_2^2,
  \end{array}
\end{equation}
which we can solve with the following theorem.
\begin{theorem}\label{thm:stoch_quad}
  Assume that $x$ is drawn from $\mathcal{N}(\mu,\Sigma)$.
  The stochastic problem~\eqref{prob:multi_step_stochastic} can be written as the deterministic problem~\eqref{prob:multi_step_eigen}, where $\bar{z}^k_j = (a_j^k)^T \mu + (a_j^k)^T \Sigma a_j^k$ and $a_j^k$ is defined as the $j$-th row of the matrix $-(I - \theta^{k-1}{\diag}(\lambda)) \cdots (I - \theta^{0}{\diag}(\lambda)) \diag \lambda^{-1} Q^T$.
  Hence, the two and three-step optimal step sizes for the stochastic problem~\eqref{prob:multi_step_stochastic} can be computed using \Thm~\ref{thm:two_step} and \Thm~\ref{thm:three_step} under the same mild conditions.
\end{theorem}
See Appendix \Sec~\ref{proof:stoch_quad} for the proof.
The upshot of \Thm~\ref{thm:stoch_quad} is that when $x$ is drawn from a known Gaussian distribution, there is no need to sample training instances; instead, we can directly solve the stochastic problem.
Moreover, the results from the previous section can directly be used.

In the following theorem, we prove a convergence rate that arises from repeatedly solving the $B$-step progressive training problem~\eqref{prob:multi_step_stochastic}. 
\begin{theorem}\label{thm:stoch_rate}
  Let $\mu$ and $L$ be the minimum and maximum eigenvalues of $P$, and let $B \in \{1,2,3\}$.
  Let $(\beta^0, \dots, \beta^{B-1})$ be Young's Chebyshev step sizes~\citep{Young1953}, and let $(\theta^0, \dots, \theta^{Bk-1})$ be learned with our optimal $B$-step progressive training to solve problem~\eqref{prob:multi_step_stochastic}.
  Then the following rate holds true:
  \begin{equation*}
    \underset{x\sim\mathcal{X}}{\mathbf{E}} \|z^{Bk}_\theta(x) - z^\star(x)\|_2 \leq \left( \max_{\mu \leq \lambda \leq L} |\Pi_{l=0}^{B-1}(1 - \beta^l \lambda)| \right)^k \underset{x\sim\mathcal{X}}{\mathbf{E}} \|z^0_\theta(x) - z^\star(x)\|_2.
  \end{equation*}
\end{theorem}
See Appendix \Sec~\ref{proof:stoch_rate} for the proof.
\Thm~\ref{thm:stoch_rate} gives a rate of convergence in expectation for our learned step sizes with progressive training with $B \in \{1,2,3\}$.
Specifically, it shows that our rate of convergence in expectation achieved through our $B$-step progressive training is at least as good as a periodic sequence of step sizes (with a period of $B$) given by Young's Chebyshev step sizes.
Note that in \Thm~\ref{thm:stoch_rate}, the initial point $z^0_\theta(x)$ must be the zero vector.

\subsection{Using gradient-based methods to solve the training problem}\label{subsec:grad_method}
Unless the algorithm of interest is gradient descent, it is challenging to solve the progressive training problems in closed-form.
Instead, we use gradient-based methods to optimize the algorithm hyperparameters.
We rely on automatic differentiation~\citep{autodiff} to differentiate through the fixed-point iterations.
We note that due to the inclusion of proximal and projection steps, there are non-differentiable mappings in our computational graph.
At these points we use subgradients to estimate directional derivatives of the loss.

It is important to restrict the allowable values that some of the hyperparameters can take for two main reasons.
First, we restrict the steady-state hyperparameters so that we can guarantee convergence.
Table~\ref{table:fp_algorithms} provides the necessary restrictions for specific algorithms.
Second, we pick sensical values over both the step-varying phase and the steady-state phase, \eg, the step size in gradient descent should always be positive.
\begin{itemize}
  \item {\bf Enforcing positivity}. To enforce positivity of the weight $\theta^k_i$ (the $i$-th hyperparameter in the $k$-th step), we set $\theta^k_i = \exp(\nu^k_i)$ and freely optimize over $\nu^k_i$.
  \item {\bf Enforcing range constraints}. To enforce that the weight $\theta^k_i$ belongs in the range $(a,b)$, we set $\theta^k_i = (b - a)/ (1 + \exp(-\nu^k_i)) + a$ and freely optimize over $\nu^k_i$.
  This is a scaled and shifted version of the sigmoid function $\sigma(\theta) = 1 / (1 + \exp(-\theta))$.
\end{itemize}

\subsection{Safeguarding mechanism}\label{subsec:safeguard}
While our method is guaranteed to converge after enough iterations, it is still possible, even if unlikely, that the iterates will become too large within the step-varying phase. 
In gradient descent and proximal gradient descent, we find that our method produces step sizes in the step-varying phase that are significantly larger than the classical ranges that guarantee convergence.
This finding aligns with a growing interest in using long step sizes to accelerate gradient descent~\citep{grimmer2023accelerated,grimmer2024provably,altschuler2023acceleration_str_cvx,altschuler_nonstr}.
To handle cases where iterates become too large, we implement a safeguarding mechanism, a common strategy in the learning to optimize literature~\citep{amos_tutorial,safeguard_convex}.
When an update gives poor performance (\eg, the objective increases by a specific factor in a single iteration), we revert to a safe, fallback update (\eg, the default step size for gradient descent). 
Once the safeguard mechanism is triggered we use the safeguarded step size from this point onwards.
This detail is important for maintaining fairness in comparison to other methods since it means that we only require a single fixed-point update in each iteration apart from at most one step where the safeguard is triggered (in this step two updates are calculated).

\section{Generalization guarantees for unseen data}\label{sec:gen}
Our convergence guarantees ensure that LAH asymptotically reaches an optimal solution, but do not provide performance guarantees for a fixed number of iterations on new instances. 
Here, we address this by presenting a method to obtain generalization guarantees, assuming access to a validation set of $N^{\rm val}$ i.i.d. samples $S = \{x_i\}_{i=1}^{N^{\rm val}}$ unseen during training. 
This procedure is applied post-training with fixed weights $\theta$.

\myparagraph{The error metric and the risk}
Recall from \Sec~\ref{sec:framework} that we evaluate LAH with a performance metric $\phi(z,x)$.
A central object in our bounds is the 0--1 error function based on the metric $\phi(z,x)$, number of algorithm steps $k$, and tolerance $\epsilon$
\begin{equation*}
  e_\theta(x) = \mathbf{1} (\phi(z^k_\theta(x), x) \ge \epsilon).
\end{equation*}
This function takes a value of one if the performance metric after $k$ steps $\phi(z^k_\theta(x), x)$ exceeds a given threshold $\epsilon$, and zero otherwise.
We aim to bound the risk $\risk(\theta)$ in terms of the empirical risk $\emprisk(\theta)$ over a validation set, with the two quantities defined as
\begin{equation}\label{eq:risk_lah}
  \risk(\theta) = \underset{x \sim \mathcal{X}}{\mathbf{E}}  e_\theta(x) = \mathbf{P}(\phi(z^k_\theta(x), x) \ge \epsilon),\quad \text{and} \quad \emprisk(\theta) = \frac{1}{N^{\rm val}} \sum_{i=1}^{N^{\rm val}} e_\theta(x_i).
\end{equation}

\myparagraph{Bounding the risk}
We adapt the method from~\citet[\Sec~5]{sambharya2024data} to obtain generalization guarantees for unseen data.
\begin{theorem}[Sample convergence bound~\citep{langford_union_prior}]\label{thm:sample_conv_bound}
  Given $\delta \in (0,1)$ and a sample dataset $S$ of size $N^{\rm val}$, with probability at least $1 - \delta$ the following bound holds:
\begin{equation}\label{eq:langford_bound}
  {\rm kl} (\emprisk(\theta) \parallel  \risk(\theta) ) \leq \frac{\log (2 / \delta)}{N^{\rm val}}.
\end{equation}
\end{theorem}
Here, the notation ${\rm kl}(q~||~p)$ is the Kullback-Liebler (KL) divergence between two Bernoulli distributions with key parameters $q$ and $p$~\citep{kl_div}; that is, ${\rm kl}(q~||~p) = q \log q/p + (1 - q) \log (1 - q)/(1 - p)$.
The inequality~\eqref{eq:langford_bound} provides an implicit bound on the risk $r_\mathcal{X}(\theta)$.
To convert the implicit bound to an explicit bound, we invert the KL divergence.
The sample convergence bound from~\eqref{eq:langford_bound} implies the inequality $p \leq {\rm kl}^{-1}(q ~|~ c) = \max \{p \in [0,1] \mid {\rm kl}(q \parallel p) \leq c \}$.
To obtain the explicit bound, we maximize over $p$.
This is a convex optimization problem in the decision variable $p \in \reals$.

A similar approach can be used to obtain lower bounds on the risk that hold with high probability.
We use the same sample convergence bound, but instead of bounding the risk in terms of the empirical risk, we bound the success rate $1 - \risk$ in terms of the empirical success rate $1 - \emprisk$.
After using \Thm~\ref{thm:sample_conv_bound} with the KL inverse, we combine the upper and lower bounds on the risk via a union bound to obtain the following sandwiched inequality which holds with probability at least $1 - 2 \delta$:
\begin{equation*}
  1 - {\rm kl}^{-1} \left(1 - \emprisk ~\bigg|~ \frac{\log (2 / \delta)}{N^{\rm val}}\right)  \leq \risk \leq {\rm kl}^{-1} \left(\emprisk ~\bigg|~ \frac{\log (2 / \delta)}{N^{\rm val}}\right).
\end{equation*}

\myparagraph{Quantile guarantees on the performance metric}
We use the bounds (both upper and lower) on the risk described above to bound the performance metric quantiles.
We achieve this by using a union-bound argument as in~\cite{sambharya2024data}.
First, we fix the number of iterations $k$, the metric $\phi$ we are interested in, and the desired confidence~$\delta$.
We discretize the entire range of possible performance metric tolerances $\epsilon$ (\eg, $10^{-10}$ to $10^5$ evenly on a log scale) into $N^{\rm tol}$ values.
For each tolerance $\epsilon$, we use the above procedure to bound the risk for that tolerance.
By virtue of a union bound, all of these risk bounds hold with probability $1 - \delta N^{\rm tol}$.
Recall that the risk with a given tolerance $\epsilon$ is the probability that the metric is less than $\epsilon$ from \Eqn~\eqref{eq:risk_lah}.
Thus, the upper bound on the risk is a valid upper bound on the quantile.
The union bound allows us to take the tightest such upper bound.
We replicate this union-bound argument to get tight  lower bounds.
Both the upper and the lower bound hold with probability $1 - 2 \delta N^{\rm tol}$ by virtue of another union bound.

\section{Numerical experiments}
In this section, we show the effectiveness of our method with many different examples.
We apply our technique to learn the hyperparameters of gradient descent in \Sec~\ref{subsec:gd_num}, proximal gradient descent in \Sec~\ref{subsec:proxgd_num}, OSQP in \Sec~\ref{subsec:osqp_num}, and SCS in \Sec~\ref{subsec:scs_num}.
The code to reproduce our results is available at
\begin{equation*}
\text{\url{https://github.com/stellatogrp/learning\_algorithm\_hyperparameters}}.
\end{equation*}
In each example, we use only $10$ training samples and evaluate with $1000$ test samples.
In all of our examples except for one, we consider $50$ time-varying steps and aim to minimize the mean square error to the ground truth after $60$ steps (\ie, $H=50$ and $K=60$).
In all of our numerical examples, we implement the progressive training regime and train $10$ algorithm steps at a time (unless otherwise indicated).
We refer the reader to Appendix~\ref{sec:num_exp_details} for details on the generalization guarantees and safeguarding.

\myparagraph{Baseline comparisons}
We compare our approach against several baselines.
First, we consider different initialization techniques with default hyperparameters in the algorithm.
\begin{itemize}
  \item \textbf{Vanilla.} We initialize the fixed-point algorithm with the zero vector.
  \item \textbf{Nearest neighbor.} The nearest-neighbor warm start initializes the test problem with an optimal solution of the nearest of the training problems measured by distance in terms of its parameter $\theta \in \reals^d$.
  In all of the examples, the distribution of problem instances is wide enough so that the nearest-neighbor warm start hardly improves (if at all) over the vanilla method.
  \item \textbf{Learned warm starts (L2WS)~\citet{l2ws}.} This technique learns a neural networks that maps the problem parameter to a high-quality warm start.
  We also use the regression loss with $10$ unrolled steps and a neural network with two layers of $500$ hidden nodes each, as these choices generally yield strong results in the experiments in~\citet{l2ws}.
\end{itemize}
Second, we consider a strategy that learns the hyperparameters of algorithm steps.
\begin{itemize}
  \item \textbf{Learned metrics (LM)~\citet{metric_learning}.} 
This approach learns the metric for Douglas-Rachford splitting~\citep{dr_splitting} to solve parametric quadratic programs by training a neural network to map the problem parameter to the metric.
  We adapt the approach to predict the step sizes for gradient and proximal gradient descent, where each step size is a vector with fixed component values across iterations.
Unlike our method, which varies hyperparameters across iterations but keeps them the same across problems, this approach keeps them constant across the iterations but can vary with the problem parameter.
  We train with $30$ fixed-point steps as this is in line with the method taken in the original paper~\citep{metric_learning}.
  We use a two layer neural network with $500$ nodes to map the problem parameter to the metric.
\end{itemize}
Table~\ref{table:num_weights} lists examples from our experiments, showing the number of weights for LAH, L2WS, and LM. The LAH method requires significantly fewer weights than L2WS and LM.
As outlined in \Sec~\ref{sec:framework}, we define the performance metric~$\phi(z^k(x),x)$ to be the suboptimality for unconstrained problems and the maximum of the primal and dual residuals for constrained problems.
We compare the performance across many values of iterations~$k$ by plotting the geometric mean of the performance metric over the test instances in each example.
We then report the number of iterations for the geometric mean of the performance metric to reach a specified tolerance in tables.


\begin{table}[!h]
  \small
  \centering
\renewcommand*{\arraystretch}{1.0}
\caption{Overview of the numerical examples.
For LAH, we learn the hyperparameters for $50$ time-varying steps (except for the logistic regression example where we have $100$) and a single steady-state step.
Thus, there are $51a$ weights to learn (except for logistic regression with $101a$ weights), where $a$ is the number of hyperparameters in a specific iteration (see Table~\ref{table:fp_algorithms}).
For L2WS and LM, the number of weights is the same since they both require a neural network to map inputs of size $d$ to outputs of size $\fplen$ (and we use the same neural network architecture).
We also report the number of weights for L2WS and LM, when the neural network does not have any hidden layers, \ie, it is a linear mapping from the parameter to the variables. Here, the number of weights is $d \fplen$.
}
\label{table:num_weights}
\adjustbox{max width=\textwidth}{
\begin{tabular}{lllllll}
\toprule
~ \begin{tabular}{@{}c@{}}numerical example \\~\end{tabular} 
& \begin{tabular}{@{}c@{}}algorithm \\~\end{tabular} 
& \begin{tabular}{@{}c@{}}parameter  \\size $d$ \end{tabular} 
& \begin{tabular}{@{}c@{}}fixed-point \\variables $\fplen$  ~\end{tabular} 
& \begin{tabular}{@{}c@{}}LAH\\num weights \end{tabular} 
& \begin{tabular}{@{}c@{}}L2WS/LM\\num weights \end{tabular} 
& \begin{tabular}{@{}c@{}}linear L2WS/LM\\num weights \end{tabular} \\
\midrule
\begin{filecontents*}{data/num_weights.csv}
problem;algo;d;n;lasco;l2ws;lm
ridge regression;gradient descent;500;1,000;51;1,001,500;500,000
logistic regression;gradient descent;157,000;785;101;79,144,285;123,245,000
lasso;proximal gradient descent;1,000;2,000;51;1,751,500;2,000,000
image deblurring;OSQP (ADMM);784;2,904;204;2,095,500;2,276,736
robust Kalman filtering;SCS (ADMM);100;1,150;255;826,500;115,000
max cut;SCS (ADMM);2,485;5,040;255;4,014,000;12,524,400
\end{filecontents*}
\csvreader[head to column names, /csv/separator=semicolon, late after line=\\]{data/num_weights.csv}{
problem=\colA,
algo=\colE,
d=\colF,
n=\colG,
lasco=\colB,
l2ws=\colC,
lm=\colD
}{\colA & \colE & \colF & \colG & \colB & \colC & \colD}
\bottomrule
\end{tabular}
}
\end{table}

\subsection{Gradient descent}\label{subsec:gd_num}
In this section, we apply our learning framework to the gradient descent algorithm from Table~\ref{table:fp_algorithms} to solve a family of ridge regression problems in \Sec~\ref{subsubsec:ridge} and a family of logistic regression problems in \Sec~\ref{subsubsec:logistic}.

\subsubsection{Ridge regression}\label{subsubsec:ridge}
In our first example, we consider the ridge regression problem
\begin{equation*}
  \begin{array}{ll}
  \label{prob:ridge}
  \mbox{minimize} & (1 / 2)\|Az - b\|_2^2 + \lambda \|z\|_2^2,
  \end{array}
\end{equation*}
where $A \in \reals^{m \times \fplen}$, $b \in \reals^m$, and $\lambda \in \reals_{++}$ are problem data and $z \in \reals^\fplen$ is the decision variable.
In this example, the parameter is $x=b$, so $A$ and $\lambda$ are the same for all problem instances.
This problem has a closed-form solution $z^\star(x) = (A^T A + \lambda I_n)^{-1} A^T x$.
Recall that~$L$ and $\mu$ are the smoothness and strong convexity parameters of the objective, and in this case are the maximum and minimum eigenvalues of $A$.
The condition number is $\kappa = L/\mu$.
This problem has a strongly convex objective function which means that gradient descent with a step size in the range $(0, 2 / L)$ linearly converges~\citep{mon_primer}.
Since this problem can be written in the form of problem~\eqref{prob:multi_step_stochastic} and the parameter is drawn from a known Gaussian, we use the results from \Sec~\ref{subsec:stoch} and directly solve the stochastic test optimization problem.
We solve the test problem \emph{exactly} with progressive training $1$, $2$, and $3$ steps at a time, and use gradient-based methods to progressively train $10$ steps at a time.

\myparagraph{Numerical example}
We set the dimensions to be $m = 500$ and $\fplen=1000$. 
We first sample the matrix $A$ with i.i.d. entries from the Normal distribution $\mathcal{N}(0, 1 / m)$.
Then we normalize $A$ so that each column of $A$ has a Euclidean norm of one.
The parameter $x$ is drawn from the Gaussian distribution $\mathcal{N}(0,I_m)$.
We take $\lambda = 0.01$.
In this example, we compare our method against data-driven approaches (the nearest neighbor, L2WS, and LM) where $N=10$ are used, and against non-data-driven methods (Nesterov's acceleration, the conjugate gradient method, and the silver step size schedule).

\myparagraph{Results}
The first baseline is gradient descent with a step size $2 / (\mu + L)$.
In this example, we compare our method against three additional non-data-driven methods.


The first is Nesterov's accelerated gradient method~\citep{nesterov} with iterates
\begin{equation*}
  y^{k+1} = z^k - \frac{4}{3L + \mu} \nabla f(z^k, x), \quad
  z^{k+1} = y^{k+1} + \biggl(\frac{ \sqrt{3 L / \mu + 1} -2}{ \sqrt{3 L / \mu + 1} + 2}\biggr)(y^{k+1} - y^k).
\end{equation*}
Nesterov's method picks the two sets of initial points to be equal, \ie, $y^0 = z^0$.

The second method is the silver step size rule for strongly convex, unconstrained minimization~\citep{altschuler2023acceleration_str_cvx}. 
The silver step size schedule where the condition number is $\kappa = L / \mu$ is constructed in the following way.
First, the sequences $u_k$ and $v_k$, initialized to $u_1 = v_1 = 1 / \kappa$, are constructed recursively with $u_k = v_{k / 2} / (\xi + \sqrt{1 + \xi^2})$ and $v_k = v_{k / 2} (\xi + \sqrt{1 + \xi^2})$ where $\xi = 1 - v_{k/2}$.
Then the silver step sizes are computed as $a_k = \psi(u_k)$ and $b_k = \psi(v_k)$, where $\psi(t) = (1 + \kappa t) / (1 + t)$.
Last, the silver step size schedule is defined recursively $h^{(k)} = (\tilde{h}^{(k / 2)}, a_k, \tilde{h}^{(k / 2)}, b_k)$ where $\tilde{h}^{(k / 2)}$ is the vector $h^{(k / 2)}$ without the final step.
The first few step size sequences are given by
\begin{equation*}
  h^{(1)} = (a_1), \quad
  h^{(2)} = (a_2, b_2), \quad
  h^{(4)} = (a_2, a_4, a_2, b_4), \quad
  h^{(8)} = (a_2, a_4, a_2, a_8, a_2, a_4, a_2, b_8),
\end{equation*}
and the iterates are updated as $z^{k+1} = z^k - \beta^k \nabla f(z^k, x)$, where $\beta^k$ is the $k$-th step size in the sequence $h^{(i)}$.
We refer the reader to~\citet[Section 3]{altschuler2023acceleration_str_cvx} for more details on the construction and the motivation of the fractal-like step size schedule.
The third is the conjugate gradient method~\citep{conjugate_gradient} which can be used to solve the linear system.
Problem~\eqref{prob:ridge} can be reformulated as solving the linear system $G z = h$ where $G = A^TA + 2 \lambda I_n$ and $h = A^T x$.
The iterates of the conjugate gradient method, initialized with $r^0 = h - G z^0$, $p^0 = r^0$, and $z^0 = 0$, are given by
\begin{equation*}
  \alpha^k = \frac{\|r^k\|_2^2}{(p^k)^T G p^k}, \quad
  z^{k+1} = z^k + \alpha^k p^k, \quad
  r^{k+1} = r^k - \alpha^k G p^k, \quad
  p^{k+1} = r^{k+1} + \frac{\|r^{k+1}\|_2^2}{\|r^{k}\|_2^2} p^k.
\end{equation*}
We initialize the Nesterov, silver step size, and conjugate gradient methods from the zero vector.
We adapt the learned metrics approach~\citep{metric_learning} to learn the weights of a neural network that maps the parameter to a vector of step sizes.

\begin{figure}[!h]
  \centering
  \includegraphics[width=0.75\linewidth]{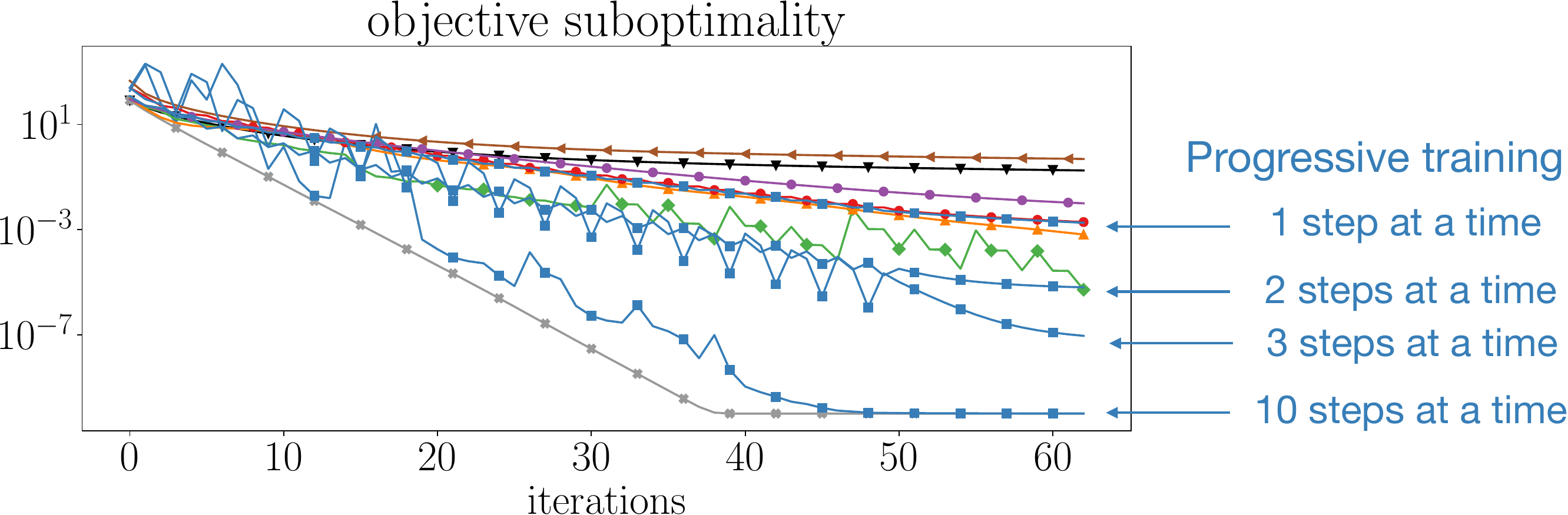}
    \\
    \legendridge
    \caption{Ridge regression results.
    The conjugate gradient method performs the best out of all of the methods.
    Our learned step sizes significantly outperform both Nesterov's method and the silver step size rule.
    Only $10$ training instances are used for each data-driven method.
    Because of this, the data-driven initialization methods, the nearest neighbor and L2WS, do not perform well.
    Considering more steps at a time in the progressive training improves the performance for LAH, but recall that the $1$, $2$, or $3$-step lookahead problems can be solved optimally, but we use gradient-based methods for the $10$-step lookahead problems.
    }
    \label{fig:ridge_results}
\end{figure}

\begin{figure}[!h]
  \centering
  \includegraphics[width=1.0\linewidth]{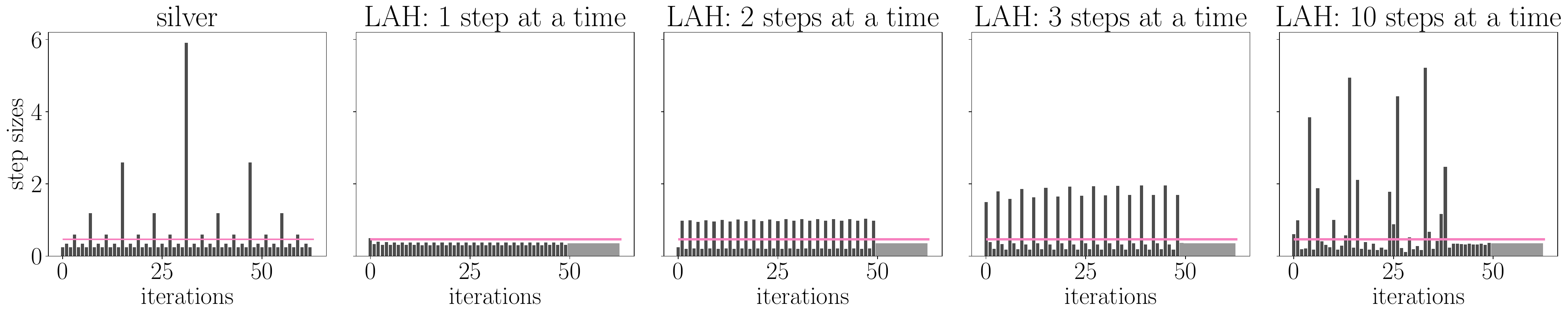}
    \\
    \legendlogisticstep
    \caption{Step sizes in gradient descent to solve the ridge regression problem.
    First on the left: silver step size schedule.
    Four on the right: our learned step sizes.
    For the first $50$ steps we learn varying step sizes in black, and for the rest, we learned a constant step size in gray.
    In pink, we show $2 / L$, the maximum constant step size that guarantees convergence.
    We observe that our learned step sizes have many short steps and several long ones -- similar to the silver step size schedule.
    }
    \label{fig:ridge_step_size_results}
\end{figure}

Figure~\ref{fig:ridge_results} and Table~\ref{tab:ridge_regression} show the behavior of our method in comparison to the baselines.
In this example, we only use $N=10$ training instances for each of the data-driven methods (apart from LAH which uses none).
The conjugate gradient method is well-known to be highly efficient for unconstrained quadratic minimization~\citep[\Sec~B.2]{acceleration_survey} and outperforms all other methods, but LAH performs the second best.
There are no quantile bounds since we directly solve the test stochastic problem.
The step sizes learned with LAH along with the silver step size schedule are depicted in Figure~\ref{fig:ridge_step_size_results}.
Interestingly, the step sizes learned with LAH exhibit several large peaks that are much larger than the maximum step size that guarantees convergence.

\begin{table}[!h]
  \centering
  \small
    \renewcommand*{\arraystretch}{1.0}
  \caption{Ridge regression results. \itersunconstrained.
    For LAH, $B$ indicates the number of steps optimized at a time. 
  }\label{tab:ridge_regression}
  \small
  \vspace*{-3mm}
  \adjustbox{max width=.95\textwidth}{
    \begin{tabular}{llllllllllll}
      \toprule
    \myCSVReaderItersRidge{./data/ridge_regression/accuracies.csv}
    \bottomrule
  \end{tabular}
  }
\end{table}

\subsubsection{Logistic regression}\label{subsubsec:logistic}
We consider the task of logistic regression used for binary classification.
Given a dataset of covariates $\{v_j\}_{j=1}^m$ where each data point $v_j \in \reals^q$ and corresponding labels $\{y_j\}_{j=1}^m$ where each label $y_j \in \{0,1\}$ this logistic regression problem can be formulated as the convex optimization problem
\begin{equation*}
  \begin{array}{ll}
  \label{prob:logisticgd_example}
  \mbox{minimize} & (1 / m) \sum_{j=1}^m y_j \log \left( \sigma(w^T v_j + b) \right) + (1 - y_j) \log \left(1 - \sigma(w^T v_j + b) \right),
  \end{array}
\end{equation*}
where the decision variables are $w \in \reals^d$ and $b \in \reals$.
Here, the mapping $\sigma : \reals \rightarrow (0, 1)$ is the sigmoid function given by $\sigma(y) = 1 / (1 + \exp(-y))$.
In this example, the parameter is the set of all of the data points and labels: $x = (v_1, \dots, v_n, y_1, \dots, y_n) \in \reals^{(q+1)m}$.
The objective function in this example is not strongly convex.

\myparagraph{Numerical example}
In our numerical experiment we consider logistic regression problems of classifying MNIST images into two classes.
To do so, for each problem, we randomly select two different classes of digits (from $0$ to $9$) and randomly select $100$ images from each class to form a dataset of $m=200$ data points.

\myparagraph{Results}
We compare our method against several additional baselines as in the previous example: Nesterov's method and the silver step size rule.
It is well-known that gradient descent with a step size of $1 / L$ achieves a convergence rate of $\mathcal{O}(1 / k)$.
Nesterov's acceleration improves upon this rate to achieve $\mathcal{O}(1 / k^2)$.
The silver step size rule achieves a convergence rate of $\mathcal{O}(1 / k^{1.2716})$, a rate in between the constant step size rate and Nesterov's accelerated rate~\citep{altschuler_nonstr}.

Nesterov's accelerated gradient method~\citep{nesterov} takes the iterates
\begin{equation*}
  y^{k+1} = z^k - \frac{1}{L} \nabla f(z^k, x), \quad
  z^{k+1} = y^{k+1} + \frac{k}{k+3}(y^{k+1} - y^k).
\end{equation*}
The two initial points are set to be equal, \ie, $y^0 = z^0$.
Given any integer $K = 2^k - 1$, the silver step size schedule of length $2K + 1$ is constructed as $h_{2K+1} = (h_K, 1 + \rho^{k-1}, h_K)$,
where $\rho = 1 + \sqrt{2}$ is the silver ratio and $h_1 = (\sqrt{2})$.
The step size in each iteration are obtained by scaling the value in the silver step size schedule with a factor of $1 / L$.
The iterates are given by $z^{k+1} = z^k - \beta^k \nabla f(z^k, x)$, where $\beta^k$ is the $k$-th step size in the sequence $h_{2K+1}$.


\begin{figure}[!h]
  \centering
  \includegraphics[width=0.75\linewidth]{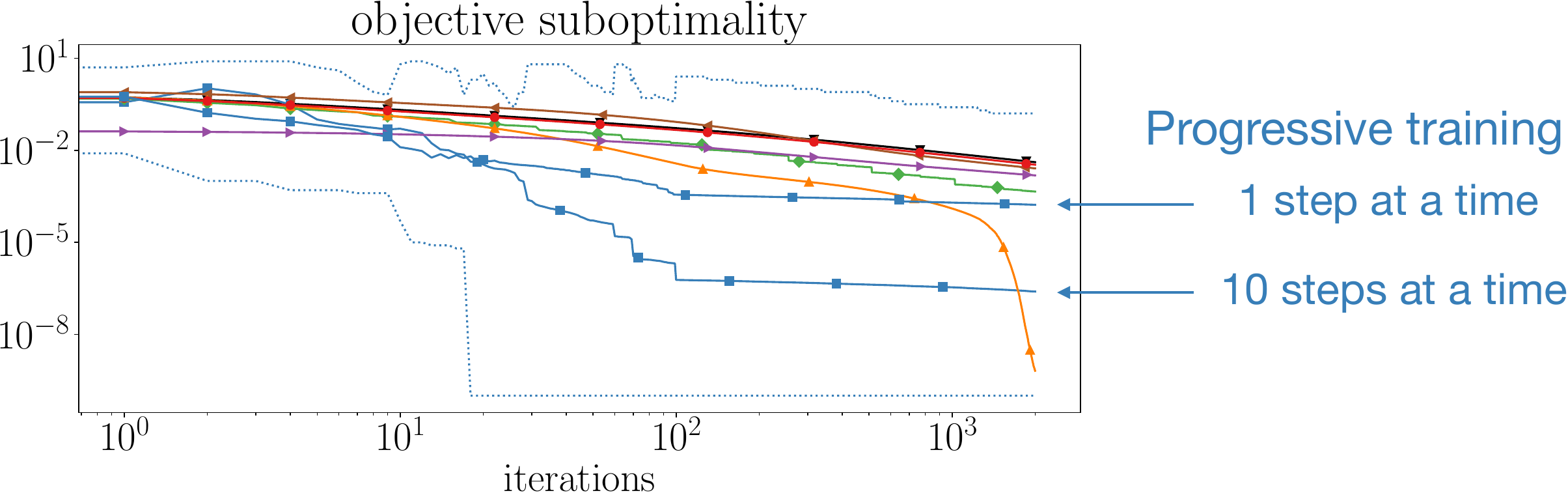}
    \\
    \legendlogistic
    \caption{Logistic regression results.
    Progressive training $10$ steps at a time reaches a geometric auboptimality average of below $10^{-5}$ within the step-varying phase (the first $100$ steps).
    In this case, we only provide the quantile bounds for progressive training $10$ steps at a time, and the bounds are wide.
    }
    \label{fig:logistic_regression_results}
\end{figure}

\begin{table}[!h]
  \centering
  \small
    \renewcommand*{\arraystretch}{1.0}
  \caption{Logistic regression results. \itersunconstrained.
  For LAH, $B$ indicates the number of steps optimized at a time.
  }\label{tab:logistic_regression}
  \small
  \vspace*{-3mm}
  \adjustbox{max width=.95\textwidth}{
    \begin{tabular}{llllllllllllll}
      \toprule
    \myCSVReaderItersLogistic{./data/logistic_regression/accuracies.csv}
    \bottomrule
  \end{tabular}
  }
\end{table}

\begin{figure}[!h]
  \centering
  \includegraphics[width=0.95\linewidth]{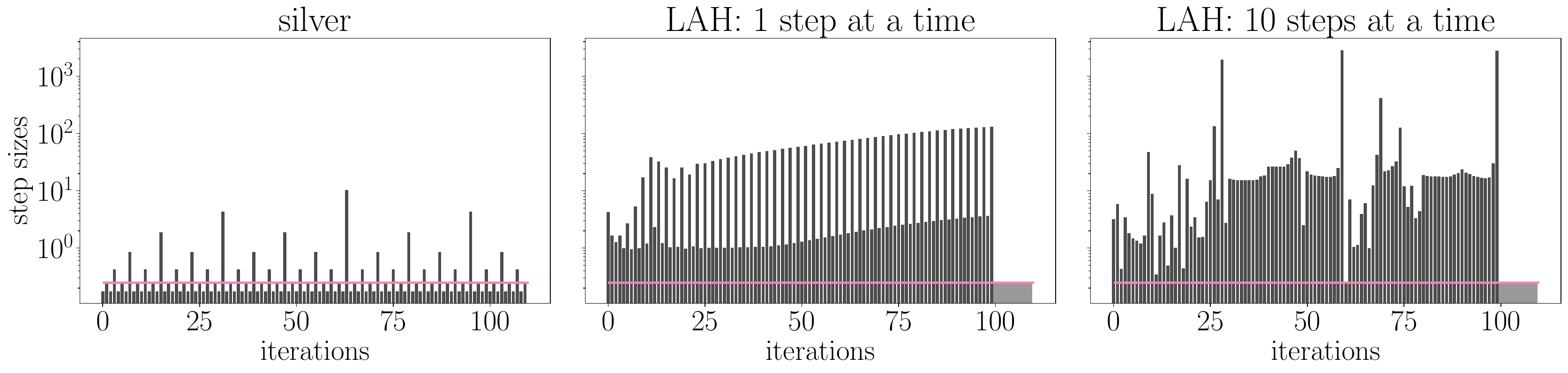}
    \\
    \legendlogisticstep
    \caption{Logistic regression step sizes.
    Many of the step sizes learned with LAH in the step-varying phase are orders of magnitude larger than the maximum constant step size that guarantees convergence ($2 / L$).
    }
    \label{fig:logistic_regression_step_sizes}
\end{figure}

\subsection{Proximal gradient descent}\label{subsec:proxgd_num}
In this section, we apply our learning framework to proximal gradient descent from Table~\ref{table:fp_algorithms} to solve a family of lasso problems in \Sec~\ref{subsubsec:lasso}.

\subsubsection{Lasso}\label{subsubsec:lasso}
In this example, we consider the lasso problem
\begin{equation}
  \begin{array}{ll}
  \label{prob:lasso}
  \mbox{minimize} & (1 / 2) \|A z - b\|_2^2 + \lambda \|z\|_1,
  \end{array}
\end{equation}
where $A \in \reals^{m \times \fplen}$, $b \in \reals^m$, and $\lambda \in \reals_{++}$ are problem data and $z \in \reals^\fplen$ is the decision variable.
In this example, the parameter is $x=b$.
Proximal gradient descent applied to problem~\eqref{prob:lasso} is the iterative shrinkage thresholding algorithm (ISTA) whose iterates are given by
\begin{equation*}
  z^{k+1} = \eta_{\lambda / L} \left(z^k - (1 / L)A^T (A z^k - x) \right),
\end{equation*}
where $\eta_\psi$ is the soft-thresholding function $\eta_\psi(z) = {\bf sign}(z) \max(0, |z| - \psi)$, and $L \in \reals_{++}$ is the largest eigenvalue of $A^T A$.
The smoothness value $L$ is the same for all problem instances.

\myparagraph{Numerical example}
To generate a family of lasso problems, we follow the setup from~\citep{l2o}.
We first sample the matrix $A$ with i.i.d. entries from the Normal distribution $\mathcal{N}(0, 1 / m)$.
Then we normalize $A$ so that each column of $A$ has a Euclidean norm of one.
To generate each problem instance, we sample a ground truth vector $z^{\rm true}$ with i.i.d. entries from the standard Normal distribution and randomly zero out $90 \%$ of the entries.
Then we set the right-hand side to be $b = A z^{\rm true} + \epsilon$ where the signal to noise ratio is set to $40$dB.
We pick the hyperparameter $\lambda$ to be $0.1$, resulting in solutions that have about $15 \%$ of their entries non-zero.
In our example, we take $m = 1000$ and $\fplen = 2000$.
This is the same ratio used in a variety of works for sparse coding problems~\citep{alista,l2o,lista_cpss}, but $m$ and $n$ are larger in our case.

\myparagraph{Results}
In this example, we also compare against the accelerated version of ISTA, known as fast ISTA (FISTA)~\citep{fista}.
The iterates of FISTA are given by
\begin{equation*}
  \small
  y^{k} = \eta_\psi(z^k - \frac{1}{L} A^T (A z^k - x)), \hspace{2mm}
  t^{k+1} = (1 / 2) \left(1 + \sqrt{1 + 4 (t^k)^2}\right), \hspace{2mm}
  z^{k+1} = y^k + \frac{t^k - 1}{t^{k+1}} (y^{k} - y^{k-1}).
\end{equation*}

\begin{figure}[!h]
  \centering
  \includegraphics[width=\figsize\linewidth]{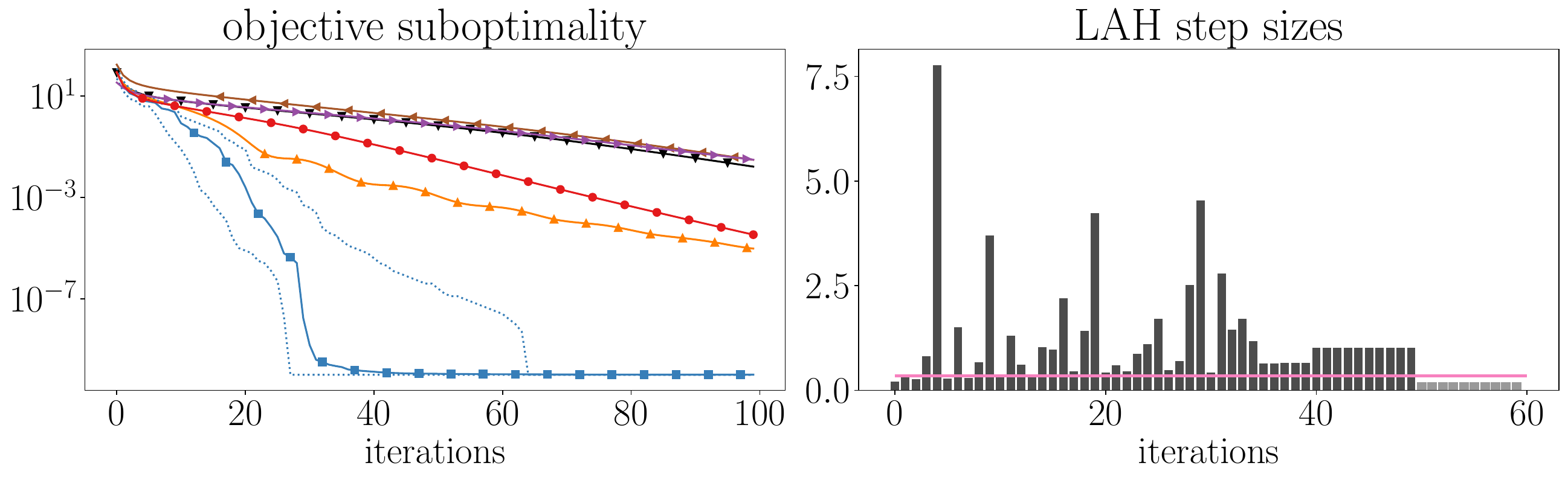}
    \legendlasso
    \caption{Left: Lasso results.
    LAH, along with its upper and lower bounds, performs the best here by a significant margin.
    On average, LAH converges to a solution within a tolerance of $10^{-8}$ within $30$ steps.
    Right: Lasso step sizes.
    The learned step size schedule has many spikes.
    LAH achieves a suboptimality lower than $10^{-9}$ after $40$ iterations, so there is little benefit in learning after this point.
    }
    \label{fig:lasso_results}
\end{figure}

Figure~\ref{fig:lasso_results} and Table~\ref{tab:lasso} show the behavior of our method in comparison to the baselines.
In this example, our method LAH performs the best by a wide margin.
The safeguarding mechanism is used in this numerical example, although it is only used in two test problems out of $1000$.
As in the case of ridge regression, the step sizes learned with LAH exhibit several large peaks that are much larger than the maximum step size that guarantees convergence.

\begin{table}[!h]
  \centering
  \small
    \renewcommand*{\arraystretch}{1.0}
  \caption{Lasso results. \itersunconstrained
  }\label{tab:lasso}
  \small
  \vspace*{-3mm}
  \adjustbox{max width=\textwidth}{
    \begin{tabular}{lllllllllll}
      \toprule
    \myCSVReaderItersLasso{./data/lasso/accuracies.csv}
    \bottomrule
  \end{tabular}
  }
\end{table}

\subsection{OSQP}\label{subsec:osqp_num}
In this section, we apply our learning framework to the OSQP algorithm~\citep{osqp} Table~\ref{table:fp_algorithms}.
We apply our method to the task of image deblurring in \Sec~\ref{subsubsec:image_deblurring}.

\subsubsection{Image deblurring}\label{subsubsec:image_deblurring}
Given a blurry image $b \in \reals^\scsnprimal$, the goal of image deblurring is to recover the original, unblurred image $y \in \reals^\scsnprimal$.
Here, the vectors $b$ and $y$ are formed by vectorizing the matrix representations of the blurry and recovered images respectively.
The image deblurring problem is
\begin{equation*}
  \begin{array}{ll}
  \label{prob:img_deblur}
  \mbox{minimize} & \|A w - b\|_2^2 + \lambda \|w\|_1 \\
  \mbox{subject to} & 0 \leq w \leq 1,
  \end{array}
\end{equation*}
where $y$ is the decision variable, $A \in \reals^{\scsnprimal \times \scsnprimal}$ is the blur operator, and $\lambda \in \reals_{++}$ is a hyperparameter that weights the quality of the image recovery the penalty that encourages sparsity of the solution.
In this example, we aim to deblur many images that are similar in nature, so the parameter is the blurry image, \ie, $x=b$.

\myparagraph{Numerical example}
We deblur EMNIST~\citep{emnist} images and follow the exact setup from~\citet{l2ws}.
We use a Gaussian blur of size $8$ and then add Gaussian noise with zero mean standard deviation $0.001$ to each pixel in an i.i.d. fashion.

\begin{figure}[!h]
  \centering
  \includegraphics[width=\figsize\linewidth]{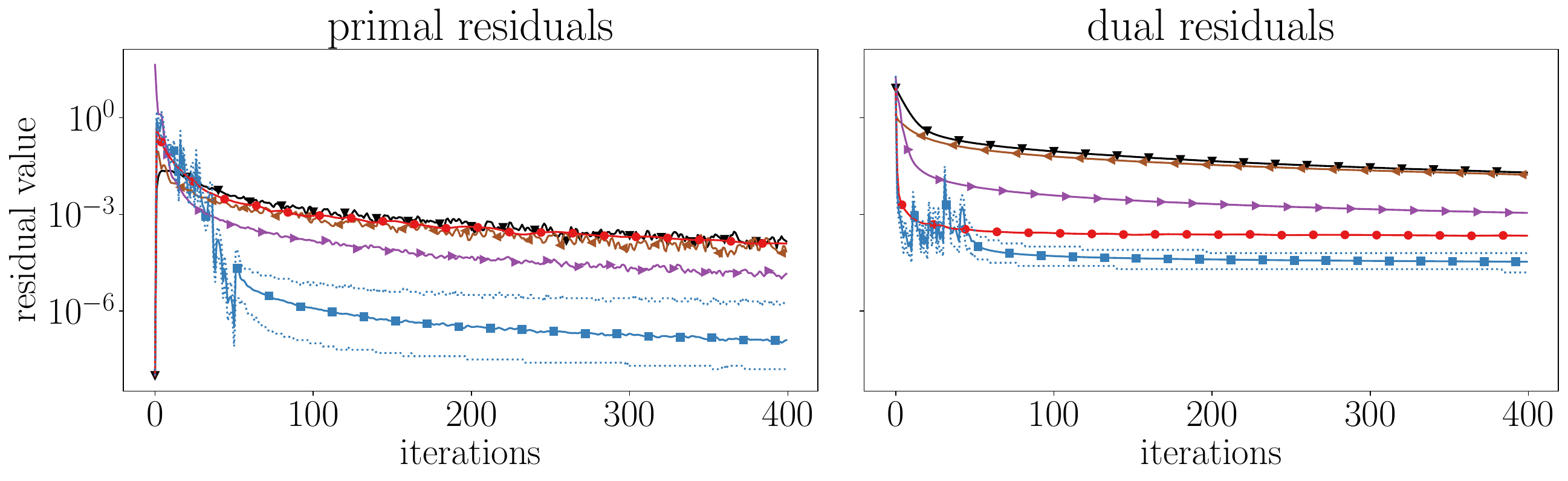}
    \\
    \legend
    \caption{Image deblurring results.
    LAH performs the best for both the primal and dual residuals.
    }
    \label{fig:mnist_results}
\end{figure}

\begin{figure}[!h]
  \centering
  \includegraphics[width=.6\linewidth]{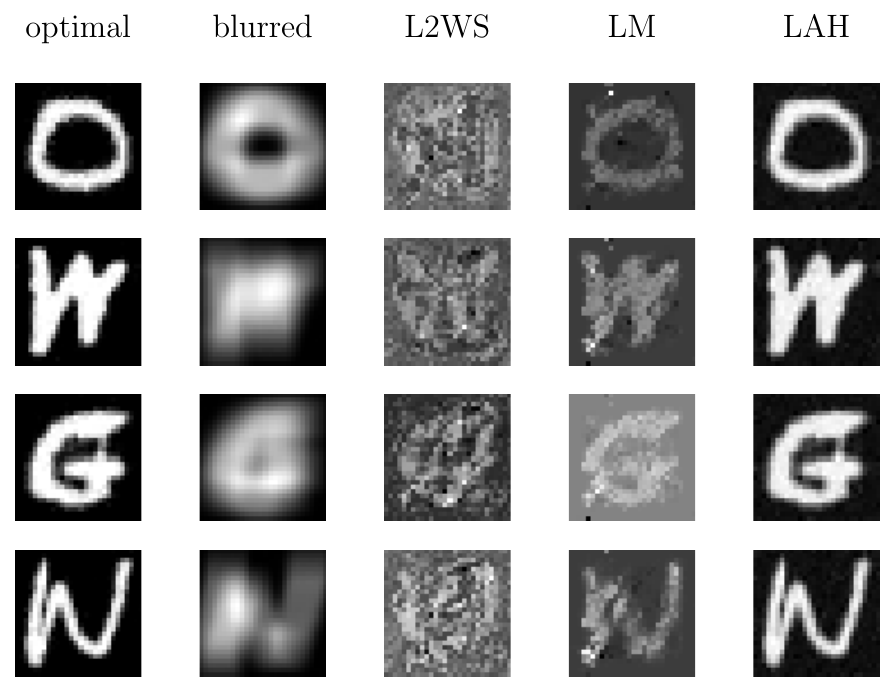}
    \\
    \caption{Image deblurring visuals.
    Each row corresponds to an unseen digit from the EMNIST dataset.
    After $3$ fixed-point steps, only LAH is able to recover the original image well.
    }
    \label{fig:mnist_visuals}
\end{figure}

\myparagraph{Results}
Figure~\ref{fig:mnist_results} and Table~\ref{tab:mnist} show the behavior of our method in comparison to the baselines.
Any method that initializes with the zero vector starts out with a primal feasible solution.
The LAH method performs the best in terms of both the primal and dual residuals using only $10$ training instances.
The L2WS method outpeforms the nearest neighbor and vanilla methods, but still takes over $6400$ iterations for both residuals to reach $0.0001$.
The LM method is able to decrease the dual residual quickly, but this comes at the expense of primal feasibility.
We visualize the effectiveness of the LAH method in Figure~\ref{fig:mnist_visuals}.

\begin{table}[!h]
  \centering
  \small
    \renewcommand*{\arraystretch}{1.0}
  \caption{Image deblurring results. \iters
  }\label{tab:mnist}
  \small
  \vspace*{-3mm}
  \adjustbox{max width=\textwidth}{
    \begin{tabular}{lllllllllllll}
    \myCSVReaderIters{./data/image_deblur/accuracies.csv}
    \bottomrule
  \end{tabular}
  }
\end{table}

\subsection{SCS}\label{subsec:scs_num}
In this section, we apply our learning framework to the SCS algorithm~\citep{scs_quadratic} from Table~\ref{table:fp_algorithms} for robust Kalman filtering in \Sec~\ref{subsubsec:rkf} and maxcut in \Sec~\ref{subsubsec:maxcut}.


\subsubsection{Robust Kalman filtering}\label{subsubsec:rkf}
Consider a linear dynamical system given by 
\begin{equation}\label{eq:rkf_state}
  s_{t+1} = A s_t + B u_t, \quad y_t = C s_t + v_t,
\end{equation}
where $s_t \in \reals^{n_s}$ is the state, $u_t \in \reals^{n_u}$ is the input vector, $y_t \in \reals^{n_o}$ is the observation, $v_t \in \reals^{n_o}$ is noise, and $A \in \reals^{n_s \times n_s}$, $B \in \reals^{n_s \times n_u}$, and $C \in \reals^{n_o \times n_s}$ are known matrices that govern the system's evolution over time.
We consider the task of Kalman filtering~\citep{kalman_filter}, \ie, estimating the state in a finite horizon in this system subject to noisy inputs $u_t$ and noise added to the measurements $v_t$.
We formulate the robust Kalman filtering problem as
\begin{equation*}
\begin{array}{ll}
\label{prob:rkf}
\mbox{minimize} & \sum_{t=1}^{T-1} \|w_t\|_2^2 + \mu \psi_{\rho} (v_t) \\
\mbox{subject to} & s_{t+1} = A s_t + B w_t \quad t=0, \dots, T-1 \\
& y_t = C s_t + v_t \quad t=0, \dots, T-1.\\
\end{array}
\end{equation*}
where the Huber penalty function~\citep{huber} parametrized by $\rho \in \reals_{++}$ that robustifies against outliers is $\psi_{\rho}(a) = \|a\|_2$ if $\|a\|_2 \leq \rho$ and $2 \rho \|a\|_2 - \rho^2$ otherwise.
The term $\mu \in \reals_{++}$ weights this penalty term.
The problem parameter is the observed noisy trajectory $x = (y_0, \dots, y_{T-1})$.

\myparagraph{Numerical example}
We use the same setup as in~\citet{l2ws}, which takes $n_s = 4$, $n_u = 2$, $n_o = 2$, $\mu = 2$, $\rho = 2$, and $T = 50$.
The dynamics matrices are given by
\begin{equation*}
\small
    \label{eq:rkf_dynamics}
    A = \begin{bmatrix}
        1 & 0 & (1 - (\gamma/2)\Delta t) \Delta t & 0\\
        0 & 1 & 0 & (1 - (\gamma/2)\Delta t) \Delta t\\
        0 & 0 & 1 - \gamma \Delta t & 0\\
        0 & 0 & 0 & 1 - \gamma \Delta t
    \end{bmatrix},\hspace{0mm}
    B = \begin{bmatrix}
        1 / 2\Delta t^2 & 0 \\
        0 & 1 / 2\Delta t^2\\
        \Delta t & 0 \\
        0 & \Delta t
    \end{bmatrix}, \hspace{0mm}
    C = \begin{bmatrix}
        1 & 0 & 0 & 0\\
        0 & 1 & 0 & 0\\
    \end{bmatrix},
\end{equation*}
where $\Delta t=0.5$ and $\gamma=0.05$ are fixed to be respectively the sampling time and the velocity dampening parameter.
We generate a family of problems by firsting generating $H$ different trajectories $N^{\rm traj}$, where $H$ is larger than the horizon length of each problem $T$.
To generate each trajectory $\{s_0^\star, \dots, s^\star_{T-1}\}$, we first let $s_0^\star = 0$ and then randomly sample the noisy inputs as $w_t \sim \mathcal{N}(0, 0.01)$ and noise in the measurements as $v_t \sim \mathcal{N}(0, 0.01)$.
The trajectories are completely defined by the dynamics equations from~\eqref{eq:rkf_state}.
Then each trajectory leads to the creation of many SOCPs.
The parameter of the first problem is $(y_0, \dots, y_{T-1})$, the second problem is $(y_1, \dots, y_{T})$, and so on.
Unlike~\citep{l2ws}, we do not take advantage of rotational invariance, making the problem instances more challenging.
In this example where problems arise sequentially, we also compare against the previous solution warm start, which initializes the current problem with the time-shifted solution of the previous problem~\citep{nonlinear_mpc}.

\begin{figure}[!h]
  \centering
  \includegraphics[width=\figsize\linewidth]{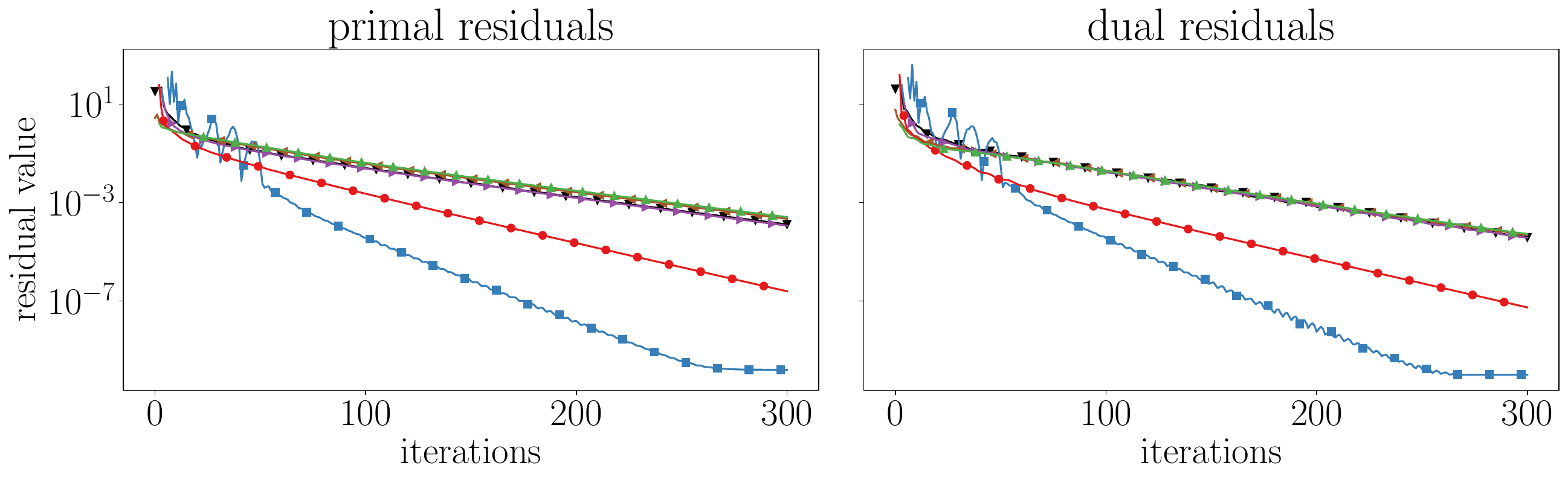}
    \\
    \legendrkf
    \caption{Robust Kalman filtering results.
    LAH performs the best, with the LM approach coming in second place.
    }
    \label{fig:robust_kalman_results}
\end{figure}

\begin{figure}[!h]
  \centering
  \includegraphics[width=.34\linewidth]{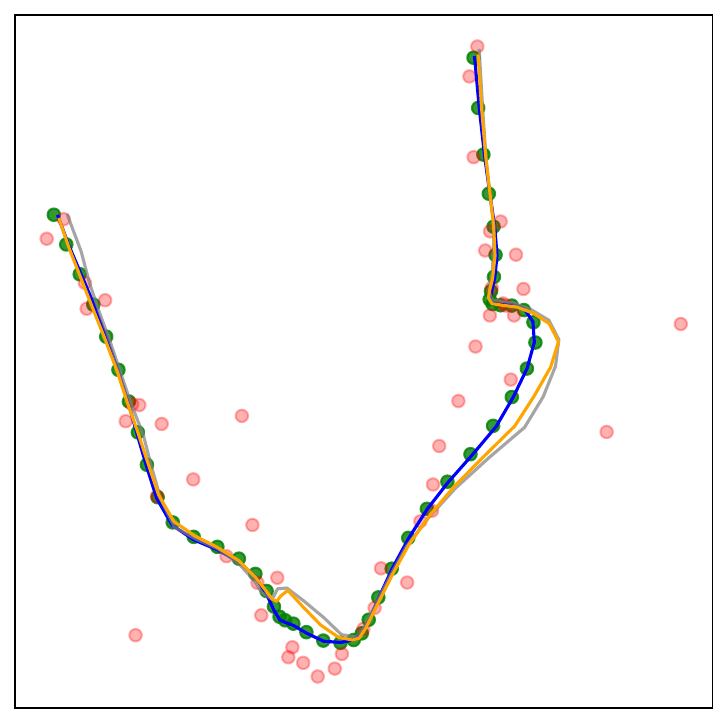}
  \includegraphics[width=.34\linewidth]{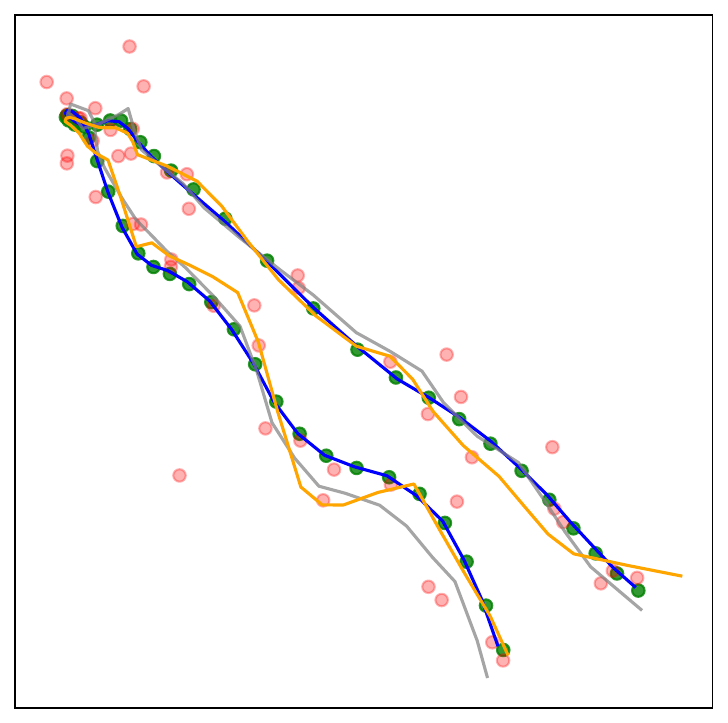}
    \rkfvisualslegend
    \caption{Robust Kalman filtering visualizations.
    The parameter for these unseen instances is the set of pink dots, and the optimal solution is the set of green dots.
    We compare the LAH, L2WS, and LM methods after $20$ fixed-point steps.
    LAH most precisely tracks the  optimal solution.
    }
    \label{fig:robust_kalman_visuals}
\end{figure}

\myparagraph{Results}
Figure~\ref{fig:robust_kalman_results} and Table~\ref{tab:rkf} show the behavior of our method in comparison to the baselines.
While the LM method with $10000$ training instances improves upon the other methods, the LAH method with only $10$ training instances improves the performance by a few orders of magnitude.
In this example, we cannot provide generalization guarantees using the result from \Sec~\ref{sec:gen} since the problem parameters are not sampled in an i.i.d. fashion.
We visualize the effectiveness of LAH in Figure~\ref{fig:robust_kalman_visuals}.

\begin{table}[htbp]
  \centering
  \small
    \renewcommand*{\arraystretch}{1.0}
  \caption{Robust Kalman filtering. \iters
  }\label{tab:rkf}
  \small
  \vspace*{-3mm}
  \adjustbox{max width=\textwidth}{
    \begin{tabular}{lllllllll}
      \toprule
    \myCSVReaderItersFull{./data/robust_kalman/accuracies.csv}
    \bottomrule
  \end{tabular}
  }
  \\
\end{table}

\subsubsection{Maxcut}\label{subsubsec:maxcut}
Consider the problem of finding the maximum cut (maxcut) in a graph (\ie, a partition of all of the nodes into two groups that maximizes the total weight across the partition). 
Given the Laplacian matrix $C \in \symm^\scsnprimal$ of the graph, the maxcut problem can be formulated as the non-convex problem of maximizing $w^T C w \quad s.t. \quad w_{i}^2 = 1 \quad i = 1, \dots, \scsnprimal$,
where $w \in \reals^\scsnprimal$ is the decision variable.
The semidefinite program
\begin{equation*}
  \begin{array}{ll}
  \label{prob:maxcut}
  \mbox{maximize} & \Tr (C W) \\
  \mbox{subject to} & W_{ii} = 1 \quad i = 1, \dots, n, \\
  & W \succeq 0,
  \end{array}
\end{equation*}
approximates the solution of this non-convex problem by lifting the vector variable $w$ to a matrix variable $W \in \symm^\scsnprimal$ and removing a rank constraint~\citep{goemans1995improved}. 
The problem parameter is the upper triangular entries of the Laplacian matrix in vectorized form, \ie, $x = {\bf triu}(C)$.

\myparagraph{Numerical example}
Our parametric family consists of random Erdos-Renyi graphs~\citep{erdds1959random} where the probability of each edge being present is $0.5$ with $70$ nodes.

\begin{figure}[!h]
  \centering
  \includegraphics[width=\figsize\linewidth]{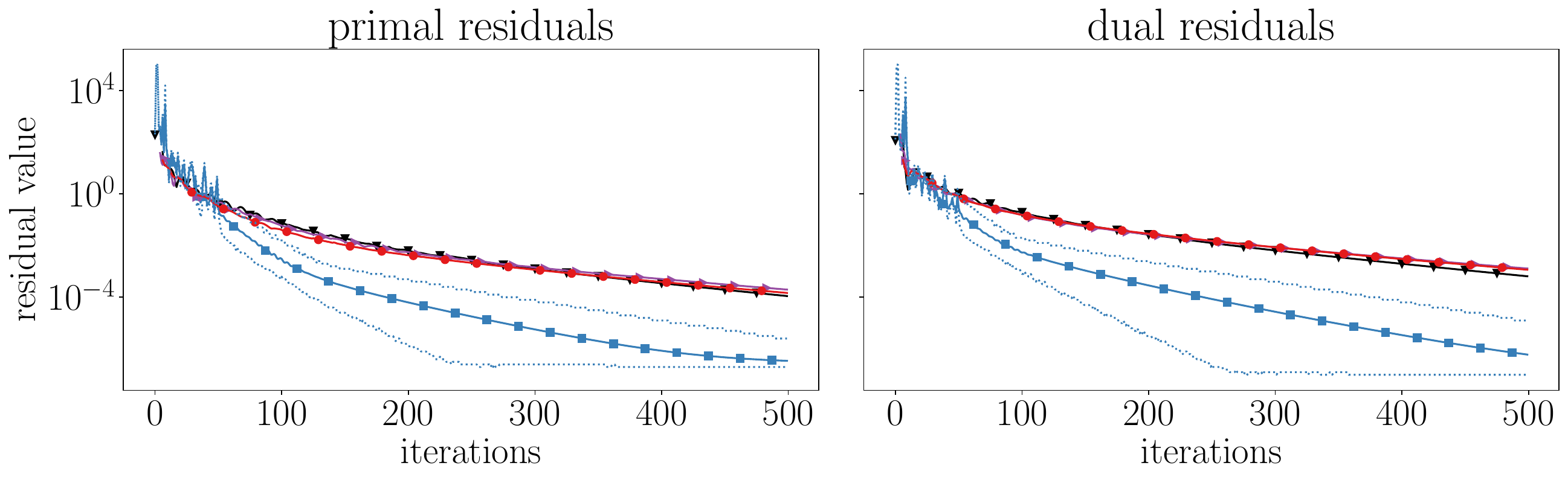}
    \\
    \legend
    \caption{Maxcut results.
    In this example, only LAH siginificantly improves the performance over the vanilla SCS method.
    \vspace*{0mm}
    }
    \label{fig:maxcut_results}
\end{figure}

\myparagraph{Results}
Figure~\ref{fig:maxcut_results} and Table~\ref{tab:maxcut} show the behavior of our method in comparison to the baselines.
LAH performs the best by a significant margin.
The $95$-th quantile upper bound guarantees also perform much better than the other methods.

\begin{table}[!h]
  \centering
  \small
    \renewcommand*{\arraystretch}{1.0}
  \caption{Maxcut. \iters
  }\label{tab:maxcut}
  \small
  \vspace*{-3mm}
  \adjustbox{max width=\textwidth}{
    \begin{tabular}{llllllllll}
      \toprule
    \myCSVReaderIters{./data/maxcut/accuracies.csv}
    \bottomrule
  \end{tabular}
  }
\end{table}

\bibliographystyle{siamplain}
\bibliography{bibliographynourl}



\ifpreprint
\appendix
\section{First-order methods}\label{sec:fom_details}
\myparagraph{Gradient descent} 
Here, $z \in \reals^n$ is the decision variable, and $f : \reals^\fplen \times \reals^d \rightarrow \reals$ is an $L$-smooth, convex objective function with respect to $z$.

\myparagraph{Proximal gradient descent} 
Here, $z \in \reals^n$ is the decision variable, $h : \reals^\fplen \times \reals^d \rightarrow \reals$ is an $L$-smooth, convex function with respect to $z$, and $g : \reals^\fplen \times \reals^d \rightarrow \reals$ is a non-smooth, convex function with respect to $z$.

\myparagraph{OSQP} 
OSQP splits the vector $\rho^k = (\rho^k_{\rm eq} \mathbf{1}_{m_{\rm eq}}, \rho^k_{\rm ineq} \mathbf{1}_{m_{\rm ineq}})$ where $m_{\rm eq}$ is the number of constraints where $l = u$, and $m_{\rm eq}$ is the number of constraints where $l < u$.
We assume that the number of equality constraints is the same for all problem instances and that the equality constraints appear first.
The primal and dual solutions at the $k$-th iteration are given by $w^k$ and $y^k = \rho(v^k - \Pi_{[l,u]}(v^k))$ respectively.
The primal residual is $\|Aw^k - \Pi_{[l,u]}(v^k)\|_2$ and the dual residual is $\|Pw^k + A^T y^k + c\|_2$.
For OSQP, the fixed-point vector $(w,v)$ is in $\reals^{\scsnprimal + m}$, \ie, $\fplen = \scsnprimal + m$.

\myparagraph{SCS} 
For SCS, $r_w \in \reals_{++}$, $r_{y_{\rm z}} \in \reals_{++}$, $r_{y_{\rm nz}} \in \reals_{++}$, and $r_\tau \in \reals_{++}$ are scaling terms that correspond to the primal variable $w$, dual variable $y$ (for both equality and inequality constraints), and the $\tau$ iterate.
We denote the number of constraints that correspond to the zero cone as $m_{\rm z}$ and the number of constraints that correspond to any other cone as $m_{\rm nz}$.
When $R = I_{m+n+1}$, we have identity-scaling.
The primal and dual solutions are given by $(w^k, y^k, s^k) = (\bar{w}^k / \tau^k, \bar{s}^k / \tau^k, \bar{y}^k / \tau^k)$.
The primal and dual residuals at the $k$-th iteration are given by $\|Aw^k + s^k - b\|_2$ and the dual residual is $\|Pw^k + A^T y^k + c\|_2$ respectively.
For SCS, the fixed-point vector $z$ is in $\reals^{\scsnprimal + m + 1}$, \ie, $\fplen = \scsnprimal + m + 1$.
We refer the reader to~\citep[\Sec~5.1]{scs_quadratic} on the details of the ${\bf root}_+$ function which involves finding the root of a quadratic equation.

\section{Numerical experiment details}\label{sec:num_exp_details}
\myparagraph{Generalization guarantees}
To generate the generalization guarantees, we use a validation set of $N^{\rm val}=1000$ samples.
We set the desired probability with $\delta=10^{-5}$ and discretize the tolerances evenly on a log scale between $10^5$ and $10^{-10}$ for a total of $N^{\rm tol}= 151$ tolerances.
Thus, with probability $1 - 2 \delta N^{\rm tol} = 0.9996$, our lower and upper quantile bounds hold simultaneously.
We report the lower bound on the $2.5$-th quantile and the upper bound on the $97.5$-th quantile.
Thus, with high probability we guarantee that at least $95 \%$ of the time the performance metric falls between the lower and the upper bound. 
We remark that we could also use our method to obtain generalization guarantees for all of the other methods, but for simplicity, only compute the bounds for our method LAH.

\myparagraph{Safeguarding}
We only use safeguarding for the logistic and lasso regression experiments to ensure that the suboptimality does not become too large during the step-varying phase.
We trigger the safeguarding mechanism if the estimated suboptimality $f(z^k(x), x) - f^\star > 10( f(z^k(x), x) - f^\star)$, where $f^\star$ is the average optimal value over the training instances.

\section{Proofs}
We first introduce Callebaut's inequality~\citep{callebaut}: for any $0 \leq s \leq t \leq 1$ and any vectors $u \in \reals^n$ and $v \in \reals^n$
\begin{equation}\label{eq:callebaut}
  \Biggl(\sum_{j=1}^n u_j^{1+s} v_j^{1-s}\Biggr) \Biggl(\sum_{j=1}^n u_j^{1-s} v_j^{1+s}\Biggr) \leq \Biggl(\sum_{j=1}^n u_j^{1+t} v_j^{1-t}\Biggr) \Biggl(\sum_{j=1}^n u_j^{1-t} v_j^{1+t}\Biggr).
\end{equation}
The more well-known Cauchy-Schwarz inequality follows by taking $s=0$ and $t=1$:
\begin{equation}\label{eq:cauchy}
  \left(\sum_{j=1}^n u_j v_j\right)^2 \leq \Biggl(\sum_{j=1}^n u_j^2\Biggr) \Biggl(\sum_{j=1}^n v_j^2\Biggr).
\end{equation}

\subsection{Proof of \Thm~\texorpdfstring{\ref{proof:pos_step}}{}}\label{sec:pos_step_proof}
For the $i$-th problem in the dataset, observe the following:
\begin{equation*}
  (\nabla f(z_\theta^k(x_i), x_i))^T (z_\theta^k(x_i) - z^\star(x_i))\geq f(z_\theta^k(x_i), x_i) - f(z^\star(x_i), x_i) \geq 0.
\end{equation*}
The first inequality uses the convexity of $f$ and the following convex property: if $f : \reals^n \rightarrow \reals$ is a convex, differentiable function then $f(y) - f(z) \geq \nabla f(z)^T (y - z)$ for all $y$ and $z$ in $\reals^n$~\citep[\Sec~3.1]{cvxbook}.
The second inequality uses the optimality of $z^\star(x_i)$.
The proof concludes by noting that the sum of non-negative values (in this case, over the $N$ problems) is non-negative, and the denominator is always positive in \Eqn~\eqref{eq:ls_sol}.
The denominator cannot be zero since there exists a $z_\theta^k(x_i)$ such that $\nabla f(z^k_\theta(x_i),x_i) \neq 0$ by assumption.

\subsection{Proof of \Thm~\texorpdfstring{\ref{thm:diff_step_sizes}}{}}\label{proof:diff_step_sizes}
To prove the theorem, we rely on a necessary second-order condition for local minima.
Without loss of generality we prove that if the first two step sizes are equal, the entire sequence of step sizes cannot be a local minimizer.
\myparagraph{Computing the partial derivatives}
We now write the partial derivatives of the training objective $r$.
The first-order derivatives are
\begin{align}\label{eq:grad_multi_step}
  \frac{\partial r(\theta^k,\dots,\theta^{k+B-1})}{\partial \theta^{k}} &= -2 \sum_{j=1}^n (1 - \theta^{k} \lambda_j) (1 - \theta^{k+1} \lambda_j)^2 \lambda_j \bar{z}_j^k \Pi_{l=k+2}^{k+B-1} (1 - \theta^{l} \lambda_j)^2  \\
  \frac{\partial r(\theta^k,\dots,\theta^{k+B-1})}{\partial \theta^{k+1}} &= -2 \sum_{j=1}^n (1 - \theta^{k} \lambda_j)^2 (1 - \theta^{k+1} \lambda_j) \lambda_j \bar{z}_j^k \Pi_{l=k+2}^{k+B-1} (1 - \theta^{l} \lambda_j)^2. 
\end{align}
The second-order derivatives are 
\begin{align}
  \frac{\partial^2 r(\theta^k,\dots,\theta^{k+B-1})}{\partial (\theta^{k})^2} &= 2 \sum_{j=1}^n (1 - \theta^{k+1} \lambda_j)^2 \lambda_j^2 \bar{z}_j^k \Pi_{l=k+2}^{k+B-1} (1 - \theta^{l} \lambda_j)^2 \\
  \frac{\partial^2 r(\theta^k,\dots,\theta^{k+B-1})}{\partial (\theta^{k+1})^2} &= 2 \sum_{j=1}^n (1 - \theta^{k} \lambda_j)^2 \lambda_j^2 \bar{z}_j^k \Pi_{l=k+2}^{k+B-1} (1 - \theta^{l} \lambda_j)^2 \\
  \frac{\partial^2 r(\theta^k,\dots,\theta^{k+B-1})}{\partial \theta^{k} \partial \theta^{k+1}} &= 4 \sum_{j=1}^n (1 - \theta^{k} \lambda_j) (1 - \theta^{k+1} \lambda_j) \lambda_j^2 \bar{z}_j^k \Pi_{l=k+2}^{k+B-1} (1 - \theta^{l} \lambda_j)^2.
\end{align}
Now by setting $\bar{\theta} = \theta^{k} = \theta^{k+1}$, we see that the upper left $2 \times 2$ entries of the Hessian are
\begin{equation*}
  2 \sum_{j=1}^n (1 - \bar{\theta} \lambda_j)^2 (1 - \theta^{k+2} \lambda_j)^2 \cdots (1 - \theta^{k+B-1} \lambda_j)^2 \lambda_j^2 \bar{z}_j^k \begin{bmatrix}
    1 & 2 \\
    2 & 1
  \end{bmatrix}.
\end{equation*}
Under the assumptions of the theorem, it is not possible for the non-negative scalar term before the $2 \times 2$ matrix to be zero.
Hence, the Hessian cannot be positive semidefinite (as Sylvester's necessary and sufficient condition to determine positive semidefiniteness is not met).
In order for a point to be a local minimum of a smooth function, it is necessary that the Hessian evaluated at that point be positive semidefinite~\citep[\Thm~2.3]{NoceWrig06}.
This completes the proof.

\subsection{Proof of \Thm~\texorpdfstring{\ref{thm:two_step}}{}}\label{proof:two_step}
First, the minimizer to problem~\eqref{prob:two_step} must exist since the objective is radially unbounded (\ie, $\lim_{\|\theta\| \to \infty} r(\theta) = \infty$).
This is a consequence of the Weierstrass Theorem ~\citep[Proposition A.8]{bertsekas1999nonlinear}.


\myparagraph{Computing the partial derivatives}
We know write the partial derivatives of the training objective $r$.
The first-order derivatives are
\begin{equation}\label{eq:grad_two_step}
  \nabla r(\theta^k, \theta^{k+1}) = 
    -2 \sum_{j=1}^n \left((1 - \theta^{k} \lambda_j) (1 - \theta^{k+1} \lambda_j)^2 \lambda_j \bar{z}_j^k,
    (1 - \theta^{k}\lambda_j)^2 (1 - \theta^{k+1} \lambda_j) \lambda_j \bar{z}_j^k \right).
\end{equation}
The second-order derivatives are 
\begin{equation}
  \nabla^2 r(\theta^k, \theta^{k+1}) = \sum_{j=1}^n \begin{bmatrix}
    -2 (1 - \theta^{k} \lambda_j) (1 - \theta^{k+1} \lambda_j)^2 \lambda_j^2 \bar{z}_j^k & 4 (1 - \theta^{k} \lambda_j) (1 - \theta^{k+1} \lambda_j) \lambda_j^2 \bar{z}_j^k\\
    4 (1 - \theta^{k} \lambda_j) (1 - \theta^{k+1} \lambda_j) \lambda_j^2 \bar{z}_j^k & -2 (1 - \theta^{k} \lambda_j)^2 (1 - \theta^{k+1} \lambda_j) \lambda_j^2 \bar{z}_j^k
  \end{bmatrix}.\!\!\!\!\!\!\!\!\!
\end{equation}
We aim to find a global minima of $r$.
Since the training objective $r$ is a smooth function, the global minima, if they exist, must be at \emph{critical points}, \ie, points $(\theta^{k},\theta^{k+1})$ where $\nabla r(\theta^{k},\theta^{k+1}) = 0$.
We break our analysis into two parts: the first where we assume that the step sizes $\theta^{k}$ and $\theta^{k+1}$ are the same, and the second where we assume they are not equal.
As proven in \Thm~\ref{thm:diff_step_sizes}, if the two step sizes are equal, that point cannot be a local minimum.
So we can skip immediately to the case of unequal step sizes.
\myparagraph{Finding the critical points in the case of unequal step sizes}
We now focus on the case where $\theta^{k} \neq \theta^{k+1}$.
We set $\nabla r(\theta^k, \theta^{k+1}) = 0$ (where $0$ is a vector in $\reals^2$) in \Eqn~\eqref{eq:grad_two_step}.
We then subtract the first of these equations from the second which yields the equation $(\theta^{k+1} - \theta^{k}) \sum_{j=1}^n (1 - \theta^{k} \lambda_j) (1 - \theta^{k+1} \lambda_j) \lambda_j^2 \bar{z}_j^k = 0$.
By assumption, $\theta^{k} \neq \theta^{k+1}$, so we can simplify this equation to
\begin{equation}\label{eq:simplify_two_step_proof}
  \sum_{j=1}^n (1 - \theta^{k} \lambda_j) (1 - \theta^{k+1} \lambda_j) \lambda_j^2 \bar{z}_j^k = 0.
\end{equation}
We now solve for $\theta^{k}$ in terms of $\theta^{k+1}$ in both \Eqn~\eqref{eq:simplify_two_step_proof} and $\partial r(\theta) / \partial \theta^k = 0$ from \Eqn~\eqref{eq:grad_two_step}.
This yields
\begin{equation*}
  \theta^{k} = \frac{\sum_{j=1}^n (1 - \theta^{k+1} \lambda_j)^2 \lambda_j \bar{z}_j^k}{\sum_{j=1}^n (1 - \theta^{k+1} \lambda_j)^2 \lambda_j^2 \bar{z}_j^k}, \quad \theta^{k} = \frac{\sum_{j=1}^n (1 - \theta^{k+1} \lambda_j) \lambda_j^2 \bar{z}_j^k}{\sum_{j=1}^n (1 - \theta^{k+1} \lambda_j) \lambda_j^3 \bar{z}_j^k}.
\end{equation*}
Setting the right hand sides of these equations equal to each other yields a cubic equation in terms of $\theta^{k+1}$.
Yet, the cubic term cancels, and we can simplify this to the following quadratic equation
\begin{equation*}
  c_2 (\theta^{k+1})^2 + c_1 \theta^{k+1}  + c_0 = 0, \quad \quad \text{where} \quad \quad c_0 = ac - b^2, \quad c_1 = bc-ad, \quad c_2 = bd - c^2.
\end{equation*}
We apply the quadratic formula to get the solutions to this equation: $\bar{\theta}^{k+1} = (-c_1 \pm \sqrt{c_1^2 - 4 c_0 c_2}) / (2c_2)$.
For a smooth function, if the Hessian evaluated at a critical point is positive definite, then that point is a strict local minimum~\citep[\Thm~2.5]{NoceWrig06}. 
Recall that our solution satisfies \Eqn~\eqref{eq:simplify_two_step_proof}, which means that the mixed second-order partial derivative is zero.
Thus the Hessian takes the form
\begin{equation*}
  \nabla^2 r(\bar{\theta}^{k}, \bar{\theta}^{k+1}) = 2 \begin{bmatrix}
    \sum_{j=1}^n (1 - \bar{\theta}^{k+1} \lambda_j)^2 \lambda_j^2 \bar{z}_j^k & 0\\
    0 & \sum_{j=1}^n (1 - \bar{\theta}^{k} \lambda_j)^2 \lambda_j^2 \bar{z}_j^k
  \end{bmatrix},
\end{equation*}
where $\bar{\theta}^k$ and $\bar{\theta}^{k+1}$ are the roots $(-c_1 \pm \sqrt{c_1^2 - 4 c_0 c_2}) / (2c_2)$ in either order.
The Hessian is positive definite since all $\lambda_j$ are positive, all $\bar{z}_j^k$ are non-negative, and by assumption of the theorem, there are at least two strictly positive $\bar{z}_j^k$.
Therefore the sufficient second-order condition is satisfied and the two solutions are local minima.
Since they are the only local minima, they are global minima.

In the rest of this proof we show that the discriminant is non-negative and that the solutions are non-negative.
By Cauchy-Schwarz~\eqref{eq:cauchy}, both $c_0$ and $c_2$ are non-negative.
The term $c_1$ is negative by the Callebaut inequality (this can be verified by taking $s=1/4$ and $t=3/4$, $u_j = \lambda_j^{9/4}\sqrt{\bar{z}^k_j}$, and $v_j = \lambda_j^{1/4}\sqrt{\bar{z}^k_j}$ and applying the Callebaut inequality from~\eqref{eq:callebaut}).
We write out the expression as follows to make the application of the Callebaut inequality clear:
\begin{align*}
  c_1 = bc-ad &=  \left(\sum_{j=1}^n \lambda_j^2 \bar{z}_j^k \right) \left( \sum_{j=1}^n \lambda_j^3 \bar{z}_j^k\right) - \left(\sum_{j=1}^n \lambda_j \bar{z}_j^k \right) \left( \sum_{j=1}^n \lambda_j^4 \bar{z}_j^k\right) \\
  &= \left(\sum_{j=1}^n \left(\lambda_j^{9/4} \sqrt{\bar{z}_j^k}\right)^{5/4}\left(\lambda_j^{1/4} \sqrt{\bar{z}_j^k} \right)^{3/4} \right) \left(\sum_{j=1}^n \left(\lambda_j^{9/4} \sqrt{\bar{z}_j^k}\right)^{3/4}\left(\lambda_j^{1/4} \sqrt{\bar{z}_j^k} \right)^{5/4} \right)\\
  &-\left(\sum_{j=1}^n \left(\lambda_j^{9/4} \sqrt{\bar{z}_j^k}\right)^{7/4}\left(\lambda_j^{1/4} \sqrt{\bar{z}_j^k} \right)^{1/4} \right) \left(\sum_{j=1}^n \left(\lambda_j^{9/4} \sqrt{\bar{z}_j^k}\right)^{1/4}\left(\lambda_j^{1/4} \sqrt{\bar{z}_j^k} \right)^{7/4} \right).
\end{align*}


\myparagraph{Non-negative discriminant}
We now show that the discriminant here is non-negative.
The discriminant is given by
\begin{align*}
  c_1^2 - 4 c_0 c_2 &= a^2d^2 + b^2c^2 + 2abcd + 4ac^3 - 4b^2c^2 - 4abcd + 4b^3d \\
  &= (ad - bc)^2 + 4ac^3  + 4 b^2 (bd - c^2) \geq 0.
\end{align*}
In the second equality, we combine like terms and factor.
In the inequality, we use the fact that $a$ and $c$ are non-negative.
We also use the fact that $bd \geq c^2$, which follows from the Cauchy-Schwarz inequality.

\myparagraph{Non-negative solutions}
It follows that the two roots to the quadratic equation are positive since $c_0$ is positive by Cauchy-Schwarz, $c_1$ is negative by the Callebaut inequality (as already shown), and $c_2$ is positive by Cauchy-Schwarz.


\subsection{Proof of \Thm~\texorpdfstring{\ref{thm:three_step}}{}}\label{proof:three_step}
As in the case of the two-step solution, we invoke the Weierstrass theorem to determine that the global minimizer to problem~\eqref{prob:three_step} exists.
We begin by writing the necessary first-order condition for $\theta^k$, $\theta^{k+1}$, and $\theta^{k+2}$:
\begin{align*}
  \sum_{j=1}^n (1 - \theta^{k+2} \lambda_j)^2 (1 - \theta^{k+1} \lambda_j)^2 (1 - \theta^{k} \lambda_j) \lambda_j \bar{z}_j^k &= 0 \\
  \sum_{j=1}^n (1 - \theta^{k+2} \lambda_j)^2 (1 - \theta^{k+1} \lambda_j) (1 - \theta^{k} \lambda_j)^2 \lambda_j \bar{z}_j^k &= 0 \\
  \sum_{j=1}^n (1 - \theta^{k+2} \lambda_j) (1 - \theta^{k+1} \lambda_j)^2 (1 - \theta^{k} \lambda_j)^2 \lambda_j \bar{z}_j^k &= 0.
\end{align*}
We add the first and second equations and then the first and third equations to get
\begin{align}
  (\theta^{k} - \theta^{k+1}) \sum_{j=1}^n (1 - \theta^{k+2} \lambda_j)^2 (1 - \theta^{k+1} \lambda_j) (1 - \theta^{k} \lambda_j) \lambda_j^2 \bar{z}_j^k &= 0, \label{eq:help_second_order_three_step} \\
  (\theta^{k} - \theta^{k+2}) \sum_{j=1}^n (1 - \theta^{k+2} \lambda_j) (1 - \theta^{k+1} \lambda_j)^2 (1 - \theta^{k} \lambda_j) \lambda_j^2 \bar{z}_j^k &= 0. \label{eq:help_second_order_three_step_two}
\end{align}
Now, we take advantage of \Thm~\ref{thm:diff_step_sizes} and use the fact that the optimal step size schedule cannot contain a pair of equal step sizes to reduce the order of the polynomials in both equations to get
\begin{align*}
  & \sum_{j=1}^n (1 - \theta^{k+2} \lambda_j)^2 (1 - \theta^{k+1} \lambda_j) (1 - \theta^{k} \lambda_j) \lambda_j^2 \bar{z}_j^k = 0 \\
  &\sum_{j=1}^n (1 - \theta^{k+2} \lambda_j) (1 - \theta^{k+1} \lambda_j)^2 (1 - \theta^{k} \lambda_j) \lambda_j^2 \bar{z}_j^k = 0.
\end{align*}
We add the above two equations and again use \Thm~\ref{thm:diff_step_sizes} to simplify and reduce the order of the polynomial to get
\begin{equation}\label{eq:three_step_cubic}
  \sum_{j=1}^n (1 - \theta^{k+2} \lambda_j) (1 - \theta^{k+1} \lambda_j) (1 - \theta^{k} \lambda_j) \lambda_j^3 \bar{z}_j^k = 0.
\end{equation}


\myparagraph{Using the two-step optimal step sizes to find the critical points}
Given a value of $\theta^k$, the step sizes $\theta^{k+1}$ and $\theta^{k+2}$ are deterministically given by \Thm~\ref{thm:two_step} in terms of $\bar{z}^{k+1}$.
Yet, by definition we have $\bar{z}_j^{k+1} = (1 - \theta^k \lambda_j)^2 \bar{z}_j^k$ for $j=1,\dots,n$, which allow us to reduce the problem to a single decision variable in $\theta^k$.
We directly write 
\begin{equation*}
  \theta^{k+2}, \theta^{k+1} = \frac{-c_1 \pm \sqrt{c_1^2 - 4 c_0 c_2}}{2c_2},
\end{equation*}
where $c_0$, $c_1$, and $c_2$ now depend on $\theta^k$:
\small
\begin{align*}
  c_0 &= \Biggl(\sum_{j=1}^n \lambda_j  (1 - \theta^k \lambda_j)^2 \bar{z}_j^k\Biggr) \Biggl(\sum_{j=1}^n \lambda_j^3 (1 - \theta^k \lambda_j)^2 \bar{z}_j^k\Biggr) - \Biggl( \sum_{j=1}^n \lambda_j^2 (1 - \theta^k \lambda_j)^2 \bar{z}_j^k\Biggr)^2 \\
  c_1 &=  \Biggl(\sum_{j=1}^n \lambda_j^2 (1 - \theta^k \lambda_j)^2 \bar{z}_j^k\Biggr) \Biggl(\sum_{j=1}^n \lambda_j^3 (1 - \theta^k \lambda_j)^2 \bar{z}_j^k\Biggr)-\Biggl(\sum_{j=1}^n \lambda_j (1 - \theta^k \lambda_j)^2 \bar{z}_j^k\Biggr) \Biggl(\sum_{j=1}^n \lambda_j^4 (1 - \theta^k \lambda_j)^2 \bar{z}_j^k\Biggr)\\
  c_2 &= -\Biggl(\sum_{j=1}^n \lambda_j^3 (1 - \theta^k \lambda_j)^2 \bar{z}_j^k \Biggr)^2 +  \Biggl(\sum_{j=1}^n \lambda_j^2 (1 - \theta^k \lambda_j)^2 \bar{z}_j^k\Biggr) \Biggl(\sum_{j=1}^n \lambda_j^4 (1 - \theta^k \lambda_j)^2 \bar{z}_j^k\Biggr).
\end{align*}
\normalsize
We define $a$, $b$, $c$, and $d$ as in \Thm~\ref{thm:two_step} and further define the following quantities:
\begin{equation*}
  e = \sum_{j=1}^n \lambda_j^5 \bar{z}_j^k, \quad f = \sum_{j=1}^n \lambda_j^6 \bar{z}_j^k.
\end{equation*}
Starting from \Eqn~\eqref{eq:three_step_cubic}, we plug in for $\theta^{k+1}$ and $\theta^{k+2}$.
After some algebraic manipulation, we get $\sum_{j=1}^n (c_2 + \lambda_j c_1 + \lambda_j^2 c_0) (1 - \theta^k \lambda_j) \lambda_j^3 \bar{z}_j^k = 0$.
At this point, we have a quintic polynomial in the variable $\theta^k$, but after further algebraic manipulation, we will find that the quintic and quartic terms cancel.
In the end, we have the following cubic equation which gives the three-step optimal step sizes:
\begin{equation*}
  d_3 (\theta^k)^3 + d_2 (\theta^k)^2 + d_1 (\theta^k) + d_0 = 0.
\end{equation*}
The coefficients are given by 
\begin{align*}
  d_0 &= ace - b^2e - ad^2 + 2bcd - c^3, \quad d_1 = c^2d - bd^2 - bce + ade + b^2f -acf \\
  d_2 &= -cd^2 + c^2e + bde - ae^2 - bcf + adf, \quad  d_3 = d^3 -  2cde + be^2 + c^2f - bdf.
\end{align*}
\myparagraph{Real roots}
The solution to the cubic equation can easily be found, for example with the Cardano method~\citep[Chapter 7]{burton2010history}.
There are two cases to the roots: either there are $3$ real roots or $1$ real root and $2$ complex roots.
We now argue that it must be the case that there are $3$ real roots.
First, we note that if $(\theta^k, \theta^{k+1}, \theta^{k+2})$ is a solution, then any permutation is also a solution.
So we can let $\theta^k$ be the real-valued solution and let $\theta^{k+1}$ and $\theta^{k+2}$ be complex-valued.
But by virtue of \Thm~\ref{thm:two_step} which we used to find the critical points, the solutions $\theta^{k+1}$ and $\theta^{k+2}$ (which depend on $\theta^k$) are real-valued.
Therefore, we reach a contradiction.

\myparagraph{Proving the critical points are optimal}
The second-order sufficient condition for a local minimum~\citep[\Thm~2.5]{NoceWrig06} can easily be verified for the critical points found.
The off-diagonal entries of the Hessian $\nabla^2 r(\theta)$ are $0$ as indicated by Equations~\eqref{eq:help_second_order_three_step} and~\eqref{eq:help_second_order_three_step_two}. 
The diagonal entries are 
\begin{equation*}
  \diag \nabla^2 r(\theta^k, \theta^{k+1}, \theta^{k+2}) = 2 \sum_{j=1}^n\begin{bmatrix}
     (1 - \theta^{k+1} \lambda_j)^2 (1 - \theta^{k+2} \lambda_j)^2 \lambda_j \bar{z}_j^k \\
    (1 - \theta^{k} \lambda_j)^2 (1 - \theta^{k+2} \lambda_j)^2 \lambda_j \bar{z}_j^k \\
    (1 - \theta^{k} \lambda_j)^2 (1 - \theta^{k+1} \lambda_j)^2 \lambda_j \bar{z}_j^k 
  \end{bmatrix}.
\end{equation*}
It is clear that the Hessian is positive definite; hence, the critical points found are local minima.
Since they are the only local minima and the global minimum is obtained, these points (there are $6$ possible permutations) are global minima.

\subsection{Proof of \Thm~\texorpdfstring{\ref{thm:stoch_quad}}{}}\label{proof:stoch_quad}
It follows from the linearity of expectation that the objective of problem~\eqref{prob:multi_step_stochastic} can be written as
\begin{equation*}
  r(\theta^k,\dots,\theta^{k+B-1}) = \sum_{j=1}^n (1 - \theta^k \lambda_j)^2 \cdots (1 - \theta^{k+B-1} \lambda_j)^2 \bar{z}_j^k,
\end{equation*}
where $\bar{z}_j^k = \underset{x \sim \mathcal{X}}{\mathbf{E}} [Q^T (z_\theta^k(x) - z^\star(x))]_j^2$.
All that remains is to simplify $\bar{z}_j^k$:
\begin{align*}
  \bar{z}^k_j &=  \underset{x \sim \mathcal{X}}{\mathbf{E}} [-Q^T Q (I - \theta^{k-1} \diag \lambda)\cdots (I - \theta^{0} \diag \lambda) \diag \lambda^{-1} Q^T x]_j^2 \\
  &= \underset{x \sim \mathcal{X}}{\mathbf{E}} [(a_j^k)^T x]^2 = \mu^T a_j^k + (a_j^k)^T \Sigma a_j^k.
\end{align*}
The first line follows from plugging in $z_\theta^0(x) = 0$ and $z^\star(x) = -P^{-1} x$.
In the second line, we use the definition of $a_j^k$ and then use the fact that $x$ is drawn from the distribution $\mathcal{N}(\mu,\Sigma)$.

\subsection{Proof of \Thm~\texorpdfstring{\ref{thm:stoch_rate}}{}}\label{proof:stoch_rate}
For every possible $x \in \reals^n$, we have
\begin{equation}\label{eq:chebyshev_rate}
  \|z_\theta^{kB}(x) - z^\star(x)\|_2^2 \leq \underset{\mu \leq \lambda \leq L}{\max} \Pi_{l=0}^{B-1} (1 - \beta^l \lambda)  \|z^{(k-1)B}_\theta(x) - z^\star(x)\|_2^2.
\end{equation}
This comes from treating $z_\theta^{k(B-1)}(x)$ as the initial points, and applying the worst-case rate of convergence for Young's Chebyshev step sizes.
We finish the proof as follows:
\begin{align*}
  \underset{x \sim \mathcal{X}}{\mathbf{E}} \|z^{kB}_{\theta}(x) - z^\star(x)\|_2^2 &\leq \underset{x \sim \mathcal{X}}{\mathbf{E}} \|(I - \beta^{k+B-1}P)\cdots(I - \beta^{k}P)(z^{(k-1)B}_{\theta}(x) - z^\star(x))\|_2^2 \\
  & \leq \underset{\mu \leq \lambda \leq L}{\max} \Pi_{l=0}^{B-1} (1 - \beta^l \lambda) \underset{x \sim \mathcal{X}}{\mathbf{E}} \|z^{(k-1)B}_\theta(x) - z^\star(x)\|_2^2.
\end{align*}
The first line holds by the optimality of $(\theta^k, \dots, \theta^{k+B-1})$ for the $B$-step progressive training problem starting from iterates $z^k_\theta(x)$.
The second line comes from the following fact; since \Eqn~\eqref{eq:chebyshev_rate} holds for every possible $x$, it holds in expectation over $x \sim \mathcal{X}$.
By applying the last inequality recursively, we conclude the proof.

\else
\fi


\end{document}